\documentclass{tran-l}

\usepackage{amsfonts}
\usepackage{amsmath}
\usepackage{graphicx}

\newtheorem{theorem}{Theorem}
\newtheorem{corollary}[theorem]{Corollary}
\newtheorem{lemma}[theorem]{Lemma}
\theoremstyle{definition}

\theoremstyle{remark}

\theoremstyle{remark}
\newtheorem{example}{Example}
\numberwithin{equation}{section}

\begin{document}
\title[Attraction and Stability of Nonlinear Ode's]{Attraction and Stability of Nonlinear Ode's using Continuous Piecewise Linear Approximations}

\author{Andr\'{e}s Garc\'{i}a and Osvaldo Agamennoni}

\address{Departamento de Ingenier\'{i}a El\'{e}ctrica y de Computadoras, Universidad Nacional del Sur, AV. Alem 1253, Bah\'{i}a Blanca, Buenos Aires, Argentina}

\email{agarcia@uns.edu.ar}

\keywords{Nonlinear systems, Stability, Attraction, Continuous Piecewise Linear.}

\begin{abstract}
In this paper, several results concerning attraction and asymptotic stability in the large of  nonlinear ordinary differential equations are presented. The main result is very simple to apply yielding a sufficient condition under which the equilibrium point (assuming a unique equilibrium) is attractive and also provides a variety of options among them the classical linearization and other existing results are special cases of the this main theorem in this paper including and extension of the well known Markus-Yamabe conjecture.

Several application examples   are presented in order to analyze the advantages and drawbacks of the proposed result and to compare such results with successful existing techniques for analysis available in the literature nowadays.
\end{abstract}

\maketitle

\section{Introduction}
Stability of nonlinear continuous system of Ordinary Differential Equations (ODE system) is a very important subject either in control design or pure theoretical analysis. It is always desirable to have systematic, necessary and sufficient conditions able to predict any class of stability where Lyapunov stability is the most common definition.


Among the available techniques nowadays, we can count mainly on Lyapunov functions as the most general method whereas many particular cases including all the available analysis for linear systems and local analysis also received considerable attention, see \cite{Arnold92} for a general reference in Lyapunov classical theory and \cite{Soliman02}, \cite{Maffei83} for successful applications. The well known Lyapunov theorem suffers the lack of a systematic counterpart making it impossible to determine if such a Lyapunov function even exists for general cases. For instance, the example by Krasowski in \cite{Bacciotti05}, pp. 32, \cite{Krasowski55} and \cite{Auslander64} shows  that sometimes such a continuous Lyapunov functions does not exists.


Several tools to construct Lyapunov functions were derived in the literature, for instance, the special case in  \cite{Papachri02} provides a methodology for polynomial systems with constraints, on the other hand, the paper in \cite{Hafstein04} was the first constructive converse Lyapunov method, in fact it was proved there that if an equilibrium point posses exponential stability behavior, then a certain feasible linear programming problem yields a Lyapunov function for such a system. Later, the work in \cite{Hafstein05} improved the result proving a similar theorem for asymptotic stability instead of exponential one. Moreover, the important reference \cite{Peet09} shows that the analysis of Lyapunov functions using polynomial approximations is topic of interest nowadays.


Alternative techniques to determine if a given system is stable (asymptotic or attractive) are the method conjectured by Aizerman (see \cite{Suddath64} and \cite{Bellman53}, pp. 86, Theorem $2$) and the result in \cite{Mangasarian63} where $x'.f(x)$ is required to be concave. Also the Corollary 2.43 in the book \cite{Chicone99} is an alternative to the Lyapunov's method, in fact this result only requires to check the stability of the jacobian of a given nonlinear autonomous system at the origin (equilibrium point), the theorem also proves that if there exists a negative upper bound for all the real part of the eigenvalues of the Jacobian matrix, then there exists a region surrounding such an upper bound where the system presents exponential behavior.\footnote{It turns out that the drawback is to determine such a region.}


Particular cases of ODE's can be successfully treated specially for low dimensions; namely $\Re^{2}$ and $\Re^{3}$ or even discontinuous systems with linear components which can be classified under the names of Hybrid Linear, Piecewise Linear, Switching Linear, etc. For these cases, several techniques have been developed; such as,  Common Quadratic Lyapunov Function, Multiple Lyapunov Functions, Poinar\'{e} Maps, Lie Algebra among many others (see for instance \cite{Johannson98},\cite{Johanson04}, \cite{DeCarlo00}, \cite{Liberzon99}, \cite{Li07} y \cite{Lin05}).


Besides the interest in stability or attraction, and as it is well known, it is neither  simple nor systematic to determine a region where any initial condition drives the system to an stable equilibrium point (domain of attraction). The most general analytic framework to find the domain of attraction requires to solve a Halmiton-Jacobi type partial differential equation (see \cite{Bacciotti05}), in this regard the paper in \cite{Ratschan06} provides such a region by defining Lyapunov (like functions which allows to prove that the trajectories never leaves the region), whereas the classical Lyapunov functions only prove the existence of such a region without providing it. Finally notice that according the section \textit{Related Work} in that paper, only two methodologies are systematic nowadays: \cite{Papachri02} and \cite{Parrilo03}.

In particular the method presented in \cite{Burhcardt07} improves the early algorithm in \cite{Ratschan06} removing some computational drawbacks, however, the algorithm starts with a given Lyapunov function, otherwise, if such a function is not provided then a linearization at the origin along with a quadratic form is used to run the method, it turns out that for certain cases (called \textit{critical cases}) the linearization exhibits zero or pure complex conjugate eigenvalues making the construction of the quadratic Lyapunov function more difficult with the necessity of some extra machinery \cite{Jyun-Horng93}.


The present paper presents two type of sufficient conditions, namely, a sufficient condition for asymptotic stability in nonlinear autonomous ODE's which makes uses of a constant negative definite matrix to produce two outcomes:

\begin{itemize}
\item{A quadratic Lyapuonv function for the nonlinear ODE being analyzed.}
\item{A subset of the domain of attraction.}
\end{itemize}

the other sufficient condition is a corollary to the previous main result, in fact instead of using a constant stable matrix, it uses the jacobian of the nonlinear vector field as a function of the state-space variables, in this case the only outcome is a subset of the domain of attraction.

The results in this paper are built on the basis of the \textit{error bounds} derived in \cite{Garcia06} and \cite{Garcia08}, where  given a  Nonlinear autonomous ODE, then an approximation of that vector field using a Continuous Piecewise Linear (CPWL) basis is considered providing dynamic bounds for an error defined as the difference between the trajectories of the Nonlinear ODE and its CPWL approximation.

This dynamic bounds (error bounds) posses the same dynamic as the CPWL approximation ODE, in this way, stability issues for the CPWL ODE are also shared by the error bounds. Realize that once the stability for the CPWL ODE is proved, then stability for the dynamic bounds is also ensured then the trajectories of the Nonlinear ODE tends to a constant or in other words also shares stability issues.

It is worth to mention the works \cite{Delanoue} and \cite{Delanoue06}, where interval analysis it is used in order to derive a method, both for asymptotic stability and the determination of some subset of the attraction domain. In fact, in \cite{Delanoue06}, an algorithm is provided to check asymptotic stability for a nonlinear ODE using successive linearization around the equilibrium at the origin. However, the stability proof is based on a quadratic Lyapunov function using the linearization aforementioned, which is not ensuring the algorithm ends in finite time for general nonlinear cases.

The work \cite{Christensen07} is very interesting and a new result, in fact, the authors of that paper derive necessary and sufficient conditions under which a given nonlinear ODE is asymptotically stable. In this way, the results in the present paper can be compared to the results in \cite{Christensen07} to notice that the present paper posses several results more general and also including the ones in that paper, however the tools in \cite{Christensen07} can be utilized together with the results in the present paper, then allowing a more extensive analysis for nonlinear systems.


This paper is organized as follows: Section \ref{The problem considered}  introduces formally the problem considered; Section \ref{Sufficient Condition in CWPL ODEs} proves a sufficient stability condition for CPWL ODE's while Section \ref{Sufficient Conditions for ODE's} provides a sufficient condition for stability or attraction in a nonlinear ODE using CPWL theory. Finally,  Section \ref{Numerical Examples} shows some examples of application.

\section{The Problem Considered}\label{The problem considered}

The goal of this paper is to derive sufficient conditions for attraction and asymptotic stability as well as quadratic Lyapunov functions for nonlinear ODE's: $\dot{x}(t)=f(x)$. The tool for such an achievement is an extension of the error bounds in \cite{Garcia06} with the formalism in \cite{Garcia08} and \cite{Garcia09}.

Thus, according to \cite{Garcia06}, if a given ODE $\dot{x(t)}=f(x)$ produces a CPWL approximation, such that this approximation posses an asymptotically stable equilibrium point, then the error bounds in the mentioned reference, also posses an asymptotically stable equilibrium point, showing that the properties of stability are shared by the given ODE and its approximation.

A correction term in the development in \cite{Garcia06} is added in this paper to take into account the possibility for the trajectories of $\dot{x}(t)=f(x)$ and its CPWL approximation to run different simplices during some interval of time.

In this way, the scenario is as follows:

\begin{gather*}
\begin{cases}
\dot{x(t)}=f(x)\quad \text{Given ODE}\\
\dot{x}_{CPWL}(t)=A^{(k)} \cdot x_{CPWL}+B^{(k)} \quad \text{CPWL Approx.}\\
x,x_{CPWL} \in \Re^{n}
\end{cases}
\end{gather*}

Recalling that the goal of the paper is to derive stability conditions of nonlinear ODE's using its CPWL approximation in the limiting caso of the grid size tending to zero\footnote{See \cite{Julian99} for a definition of grid size.}, this is only possible if the trajectories of both: nonlinear and CPWL are \textit{close each other}. To make precise what close each other means, the concept of error bounds used in \cite{Garcia06} and \cite{Garcia09} will be utilized.

In this way, the formalism in \cite{Garcia09} yields the following error bounds:

\begin{equation}\label{Error Bound Dynamics}
\begin{cases}
\dot{x(t)}=f(x)\quad \text{Given ODE}\\
\dot{x}_{CPWL}(t)=A^{(k)} \cdot x_{CPWL}+B^{(k)} \quad \text{CPWL Approx.}\\
min_{\{1,2 \}} \left\{E^{*1}(t), E^{*2}(t)\right\}\leq E^{(k)}(t) \leq max_{\{1,2 \}} \left\{E^{*1}(t), E^{*2}(t)\right\}\\
\dot{E}^{*1}(t)=A^{(i)} \cdot E^{*1}+\xi^{ik}(t)-\lambda\\
\dot{E}^{*2}(t)=A^{(i)} \cdot E^{*2}+\xi^{ik}(t)+\lambda\\
j=1,..,n
\end{cases}
\end{equation}

where $\{E^{*1},E^{*2}\}\in \Re^{n\times 1},\quad \lambda \in \Re^{n\times 1}$, $\xi^{(ik)}=\left( A^{(k)}-A^{(i)} \right)\cdot x_{CPWL}+\left( B^{(k)}-B^{(i)} \right)$ and $E^{(k)}(t)=x(t)-x_{CPWL}(t)$ is the error vector.

Notice that extra term $\xi^{(ik)}(t)$ takes into account for the possibility of the trajectories $x(t)$ running the simplex $i^{th}$ and $x_{CPWL}$ running the simplex $k^{th}$.

On the other hand, to reach the goal in this paper, it is necessary to make sure that bounds in (\ref{Error Bound Dynamics}): $min_{\{1,2 \}} \left\{E^{*1}(t), E^{*2}(t)\right\}$ and $max_{\{1,2 \}} \left\{E^{*1}(t), E^{*2}(t)\right\}$ do not grows to infinity. Clearly, if the dynamics of such a bounds: $\dot{E}^{*1,2}(t)=A^{(i)} \cdot E^{*1,2}+\xi^{ik}(t)\pm \lambda$ are attractive to one equilibrium point, then in each simplex the stability issues of the CPWL ODE is shared by the nonlinear one.

Since the error bounds in (\ref{Error Bound Dynamics}) can be rewritten in matrix block form:

\begin{equation}\label{Error Bounds Dynamics Block Form}
\underbrace{\begin{bmatrix}
\dot{E}^{*1,2}(t)\\
\dot{x}^{(k)}(t)
\end{bmatrix}}_{\dot{X}^{(k)}}=
\underbrace{\begin{bmatrix}
A^{(i)} & A^{(k)}-A{(i)}\\
0 &  A^{(k)}
\end{bmatrix}}_{\bar{A}^{(k)}} \cdot X^{(k)}+
\begin{bmatrix}
B^{(k)}-B^{(i)}\\
B^{(k)}
\end{bmatrix}+
\begin{bmatrix}
\pm \lambda\\
0
\end{bmatrix}
\end{equation}

The theorem in next subsection asserts that a CPWL ODE is attractive to one of its equilibrium points once all the matrices $\bar{A}^{(k)}$ are definite negative. In this context, (\ref{Error Bounds Dynamics Block Form}) will conduct to the desired goal in this paper in the view of the well known Theorem for block matrices: the eigenvalues of matrix $\bar{A}^{(k)}$ is the junction of the eigenvalues of matrices $\{A^{(k)},A^{(i)} \}$.

\subsection{A Sufficient Condition for Asymptotic stability in CPWL ODE's}\label{Sufficient Condition in CWPL ODEs}

In this section a theorem concerning attraction for CPWL ODE's is recalled. In fact, the theorem was proved in \cite{Garcia08} and it is important because it provides a tool which will be invariant when the grid size of the CPWL approximation simplicial partition tend to zero (as it will be shown in section \ref{Sufficient Conditions for ODE's}).

\begin{theorem}\label{Stability Criterion}
Given a CPWL ODE: $\quad \dot{x}=A^{(i)} \cdot x^{(i)} + B^{(i)}$ defined in an set $\Omega$ which posses no equilibrium points over the frontiers\footnote{This means that the equilibriums of each individual linear system in each simplex are outside the borders.}. If, moreover, for each simplex $i^{th}$, $A^{(i)}$ is negative definite and $\Omega$ is an invariant set\footnote{See \cite{Michel08}, pp. 73 and pp. 76 for a further reading.}, then the CPWL ODE is attractive to at least one of its equilibrium points wiht the initial conditions: $x(0)\in\Omega$.
\end{theorem}

Whit this Theorem and considering the grid size of a CPWL ODE approximating a given nonlinear ODE, tending to zero (the CPWL approximation tend to the nonlinear vector field, \cite{Julian99}) , it is going to be proved a theorem providing sufficient conditions for attraction and stability in a systematic way.


\section{The Main Result: A Sufficient Condition for Stability of NonLinear ODE's}\label{Sufficient Conditions for ODE's}

Up to now,  conditions for a CPWL ODE being asymptotically stable were derived, in particular when these ODE's come from a CPWL approximation of a nonlinear vector field, the quality of the approximation improves as the grid size tend to zero (see \cite{Julian99}).

In order to apply this results to predict stability/attraction in Nonlinear ODE's, two main research streams can be followed\footnote{From now on \textit{stability} will mean stability in the sense of Lyapunov (see \cite{Michel08}, pp. 76-82 for definitions of Lyapunov stability).}:

\begin{itemize}
\item{Given an ODE: $\dot{x(t)}=f(x)$ and a CPWL approximation of it, analyzing the stability of each matrix $A^{(i)}$ (this is the research line presented in \cite{Johanson04}).}
\item{Given an ODE: $\dot{x(t)}=f(x)$, utilizing the CPWL theoretical tools (like the ones in \cite{Julian99}), to derive conditions under which there exists a CPWL approximation (without explicit calculation of it), such that each matrix $A^{(i)}$ is Hurwitz, then ensuring stability or attraction in the view of Theorem \ref{Stability Criterion}.}
\end{itemize}

Following the second research line above, it is clear that the sufficient conditions in Theorem \ref{Stability Criterion}, require the following:

\begin{itemize}
\item{To derive conditions using $f(x)$ such that there exists a CPWL approximation for certain degree of accuracy where each matrix $A^{(i)}$ is definite negative.}
\item{To find an invariant set $\Omega$ for $\dot{x(t)}=f(x)$}
\end{itemize}

In order to provide a smooth introduction for the main result in this paper it is needed some extra machinery. In this way, the first condition above, requires the following Lemma:

\begin{lemma}[Stability Conditions for Linear Systems]\label{Stability of A}
Given a Linear ODE  $\dot{x(t)}=A\cdot x +B$ such that:

\begin{equation*}
\begin{cases}
\mid -F_{jk}+A_{jk} +B_{j} \mid \leq \bar{\lambda}^{*}\\
\mid B_{j} \mid \leq \widetilde{\lambda}^{*}\\
j,k=1,2,\ldots,n\\
F'\cdot P+P \cdot F=-I\\
R=
\begin{bmatrix}
\mid P \mid \cdot(\bar{\lambda}^{*}+\widetilde{\lambda}^{*})+(\bar{\lambda}^{*}+\widetilde{\lambda}^{*})' \cdot \mid P \mid
\end{bmatrix}\\
\lambda_{R} <1 \quad \text{The biggest eigenvalue of $R$}\\
\text{$\mid \cdot \mid$ means the absolute value of each element}\\
\end{cases}
\end{equation*}

Then, $A$ is negative definite-see \cite{Bellman65} for details on negative definiteness.
\end{lemma}

\begin{proof}
    The proof is in the Appendix.
\end{proof}

It is important to remember that a positive definite matrix $P$ and symmetric such that $F'\cdot P+P \cdot F=-I$, always exists if and only if $F$ is Hurwitz (see \cite{Kwatny00}, pp. 22).

On the other hand, an important corollary which will be useful in what follows is the fact that $V(x)=x'\cdot P \cdot x$ constitutes a Lyapunov function for $\dot{x(t)}=A\cdot x +B$:

\begin{corollary}\label{Lyapunov functions for Linear Systems}
Under the hypothesis of Lemma \ref{Stability of A}, $V(x)=x'\cdot P \cdot x$ is a Lyapunov function for $\dot{x(t)}=A\cdot x +B$.
\end{corollary}

\begin{proof}
    The proof is in the Appendix.
\end{proof}

Using these results for pure Linear systems and taking into account that all the development in this paper is supported by CPWL approximations, which are in essence a collection of linear systems, it is possible to prove the main result:

\begin{theorem}[Main Result]\label{Main Theorem}
Let  $\dot{x(t)}=f(x),\quad,\quad f(0)=0,\quad f\in \Re^{n}$ a Nonlinear ODE  such that $f(x)$ is Lipschitz continuous in a set $\Omega$. If given three matrices $\{F,\widetilde{\lambda}^{*},\bar{\lambda}^{*} \}$:

\begin{equation*}
\begin{cases}
F=[F_{1}\quad F_{2}\quad \ldots \quad F_{n}]\\
\bar{\lambda}^{*}=[\bar{\lambda}_{1}^{*}\quad \bar{\lambda}_{2}^{*}\quad \ldots \quad \bar{\lambda}_{n}^{*}]\\
\widetilde{\lambda}^{*}=[\widetilde{\lambda}_{1}^{*}\quad \widetilde{\lambda}_{2}^{*}\quad \ldots \quad \widetilde{\lambda}_{n}^{*}]\\
%
%
F'\cdot P+P \cdot F=-I\\
R=
\begin{bmatrix}
\mid P \mid \cdot(\bar{\lambda}^{*}+\widetilde{\lambda}^{*})+(\bar{\lambda}^{*}+\widetilde{\lambda}^{*})' \cdot \mid P \mid
\end{bmatrix}\\
\lambda_{R} <1 \quad \text{The biggest eigenvalue of $R$}\\
\text{$\mid \cdot \mid$ means the absolute value of each element}\\
\forall x\in \Omega\\
\Omega=\{x: \mid \frac{\partial{f(x)}}{\partial{x}}\cdot x \mid \leq min_{i=1,\ldots,n} \quad \widetilde{\lambda}_{i}^{*},\quad
\mid \frac{\partial{f(x)}}{\partial{x}} -F \mid \leq (\bar{\lambda}^{*}-\widetilde{\lambda}^{*})\\
\end{cases}
\end{equation*}

where $\lambda_{j}$ is the maximum error between $f(x)$ and the CPWL approximation in the spirit of \cite{Julian99} and $\breve{e}_{i}=[0\quad 0 \quad \ldots \underbrace{1}_{i^{th} position} 0 \ldots 0]'$.

Then, $\dot{x(t)}=f(x)$ is attractive to the origin in $\Omega$. Moreover, $\Omega$ is a domain of attraction.
\end{theorem}

\begin{proof}
    The proof is in the Appendix.
\end{proof}

This Theorem provides sufficient conditions for attraction in a region $\Omega$, so in order to extend this result to asymptotic stability, it will necessary the following Lemma:

\begin{lemma}[Lyapunov Functions]\label{Corollary of the Main Result}
Given an ODE: $\dot{x(t)}=f(x),\quad f\in \Re^{n}$, f Lipschitz continuous with $f(0)=0$. If given three matrices $\{F,\widetilde{\lambda}^{*},\bar{\lambda}^{*} \}$:

\begin{equation*}
\begin{cases}
F=[F_{1}\quad F_{2}\quad \ldots \quad F_{n}]\\
\bar{\lambda}^{*}=[\bar{\lambda}_{1}^{*}\quad \bar{\lambda}_{2}^{*}\quad \ldots \quad \bar{\lambda}_{n}^{*}]\\
\widetilde{\lambda}^{*}=[\widetilde{\lambda}_{1}^{*}\quad \widetilde{\lambda}_{2}^{*}\quad \ldots \quad \widetilde{\lambda}_{n}^{*}]\\
%
%
F'\cdot P+P \cdot F=-I\\
R=
\begin{bmatrix}
\mid P \mid \cdot(\bar{\lambda}^{*}+\widetilde{\lambda}^{*})+(\bar{\lambda}^{*}+\widetilde{\lambda}^{*})' \cdot \mid P \mid
\end{bmatrix}\\
\lambda_{R} <1 \quad \text{The biggest eigenvalue of $R$}\\
\frac{\partial{P}}{\partial{x}}=0
\end{cases}
\end{equation*}

then $V(x)=x'\cdot P \cdot x$ is a Lyapunov function for $f(x)$ for all $x \in \Omega$,with $\Omega$ as follows:

\begin{equation*}
\Omega=\{x: \mid \frac{\partial{f(x)}}{\partial{x}}\cdot x \mid \leq min_{i=1,\ldots,n} \quad \widetilde{\lambda}_{i}^{*},\quad
\mid \frac{\partial{f(x)}}{\partial{x}} -F \mid \leq (\bar{\lambda}^{*}-\widetilde{\lambda}^{*})\\
\end{equation*}
\end{lemma}

\begin{proof}
    The proof is in the Appendix.
\end{proof}

This Lemma defines a set $\Omega$ where a Lyapunov function $V(x)=x'\cdot P \cdot x$ is available. However, as it is well known, the set $\Omega$ could be a subset of the truly domain of attraction. In this way, to find the complete domain of attraction, the optimization algorithm in \cite{Burhcardt07} could be used, with the Lyapunov function $V(x)=x'\cdot P \cdot x$.

Clearly, Theorem \ref{Main Theorem} and Lemma \ref{Corollary of the Main Result} prove trivially the following result:

\begin{theorem}[Sufficient Conditions for Asymptotic Stability]\label{Asymptotic Stability}
Under the conditions on Lemma \ref{Corollary of the Main Result}, $\dot{x}(t)=f(x)$ is asymptotically stable in the set:

\begin{equation*}
  \Omega=\{ \mid \frac{\partial{f(x)}}{\partial{x}}\cdot x \mid \leq min_{i=1,\ldots,n} \quad \widetilde{\lambda}_{i}^{*} \}
\end{equation*}
\end{theorem}

In this way, the first possible choice for $F$ which always yields a non-empty region $\Omega$ containing the equilibrium at the origin, is $F=\frac{\partial{f(x)}}{\partial{x}}|_{x=0}$. This is not more than the classic linearization at the origin but this time, with a known validity region.

A second choice for $F$ is: $F=\frac{\partial{f(x)}}{\partial{x}}$, providing Sufficient and systematic conditions under which $\dot{x(t)}=f(x)$ is attractive to the origin:

\begin{corollary}[Sufficient Conditions for Attraction]\label{Systematic Conditions for Attraction}
Let $\dot{x(t)}=f(x),\quad,\quad f(0)=0,\quad f\in \Re^{n}$ a Nonlinear ODE  such that $f(x)$ is Lipschitz continuous in a region $\Omega$ given by:

\begin{equation*}
\Omega=\{x:\lambda_{R}(x)<1\}
\end{equation*}

where $R=2\cdot(\underbrace{\mid P \mid \cdot \mid \frac{\partial{f(x)}}{\partial{x}}\cdot x \mid \cdot [1\quad 1 \ldots 1]}_{H(x)}+H(x)')$ and $\lambda_{R}$ is the biggest eigenvalue of $R$ providing $P \cdot \frac{\partial{f(x)}}{\partial{x}}+\frac{\partial{f(x)}}{\partial{x}}' \cdot P=-I$. Then, $\dot{x(t)}=f(x)$ is attractive to the origin for all $x\in\Omega$.
\end{corollary}

On the other hand is interesting to compare our test of stability in Theorem \ref{Main Theorem} with Theorem 7.1 in \cite{Johanson04}, pp. 115-116, where the following result is presented.:

\begin{theorem}
Let z(t) be a piecewise $C^{1}$ trajectory of the system $\dot{x(t)}=f(x)$ and assume that:
\begin{equation*}
\| f(x) - A^{(i)}\cdot x -B^{(i)} \| \leq \varepsilon^{i} \cdot \| x \|^{2}
\end{equation*}

If there exists numbers $yi > 0$ and symmetric matrices $U, W^{i}$ with non-negative entries, and a symmetric matrix $T$ such that $Pi=F^{i'} \cdot T \cdot F^{(i)}$ and $\overline{Pi}=\overline{F}^{(i)'} \cdot T \cdot \overline{F}^{i}$
E~TT Fi satisfy:

\begin{equation*}
\begin{cases}
-2 \cdot \varepsilon^{i} \cdot \gamma^{(i)} \cdot I > A^{(i)} \cdot P^{(i)} +P^{(i)} \cdot A^{(i)}+E^{(i)'}\cdot U^{(i)} \cdot E^{(i)}\\
E^{(i)'}\cdot W^{(i)} \cdot E^{(i)}< P^{(i)} <\gamma^{(i)} \cdot I
\end{cases}
\end{equation*}

for $i$ in the first simplex and:

\begin{equation*}
\begin{cases}
-2 \cdot \varepsilon^{(i)} \cdot \gamma^{i} \cdot I > \overline{A}^{(i)} \cdot \overline{P}^{(i)} + \overline{P}P^{(i)} \cdot \overline{A}^{(i)}+\overline{E}^{(i)'}\cdot U^{(i)} \cdot \overline{E}^{(i)}\\
\overline{E}^{(i)'}\cdot W^{(i)} \cdot \overline{E}^{i}< P^{(i)} <\gamma^{(i)} \cdot I
\end{cases}
\end{equation*}

for $i$ in the next simplex, then $x(t)$ tends to zero exponentially.
\end{theorem}

Notice that in order to use this theorem and as is pointed out in \cite{Johanson04}, we need to run a CPWL approximation numerically and then extracting from that the matrices $A^{(i)}, B^{(i)}$ to conclude stability for our given nonlinear ODE, the point is how to ensure a correct simplicial division to get precise conclusions. In that sense our result is much more stronger since we only need the evaluation of the matrix F and the Lipschitz coefficients if an estimation of the domain of attraction has to be determined, without the necessity of matrices $A^{(i)}, B^{(i)}$.

We leave this section with a example of application of Corollary \ref{Systematic Conditions for Attraction} and Theorem \ref{Main Theorem} utilizing the jacobian at the origin.\footnote{More examples are given later in Section \ref{Numerical Examples}}:

\textbf{\textit{Example}}

The well known system of Hopf in normal form is:

\begin{equation*}
\begin{cases}
\dot{x} = f_{1}(x,y)=\alpha \cdot x- y + x \cdot (x^{2}+y^{2})\\
\dot{y} = f_{2}(x,y)=x + \alpha \cdot y + y \cdot (x^{2}+y^{2})
\end{cases}
\end{equation*}

Calculating the jacobian matrix:

\begin{equation*}
\frac{\partial{f(x)}}{\partial{x}}=
\begin{bmatrix}
\alpha+3\cdot x^{2}+y^{2} & -1+2\cdot x \cdot y\\
1 + 2 \cdot x \cdot y & \alpha+x^{2}+y^{2}
\end{bmatrix}
\end{equation*}

then, the Jacobian at the origin leads:

\begin{equation*}
\frac{\partial{f(x)}}{\partial{x}}=
\begin{bmatrix}
\alpha & -1\\
1 & \alpha
\end{bmatrix}
\end{equation*}

The eigenvalues are $\alpha \pm \sqrt(-1)$, which indicates the necessity for $\alpha <0$. In this way, the conditions on Theorem \ref{Main Theorem} yields:

\begin{equation}\label{Conditions Hopf I}
\begin{cases}
\mid \frac{\partial{f(x)}}{\partial{x}}-F \mid \leq \bar{\lambda}^{*}-\widetilde{\lambda}^{*} \Leftrightarrow
\left |
\begin{matrix}
3\cdot x^{2}+y^{2} & 2 \cdot x \cdot y\\
2 \cdot x \cdot y & 3\cdot y^{2}+x^{2}
\end{matrix} \right | \leq \bar{\lambda}^{*}-\widetilde{\lambda}^{*}\\
\mid \frac{\partial{f(x)}}{\partial{x}}\cdot x \mid \leq min_{i=1,2} \quad \widetilde{\lambda}^{*}_{i} \Leftrightarrow
\left |
\begin{matrix}
(\alpha \cdot x - y) + 3 \cdot x \cdot (x^{2}+y^{2})\\
(\alpha \cdot y + x) + 3 \cdot y \cdot (x^{2}+y^{2})
\end{matrix} \right | \leq min_{i=1,2} \quad \widetilde{\lambda}^{*}_{i}
\end{cases}
\end{equation}

On the other hand, if $\bar{\lambda}^{*},\widetilde{\lambda}^{*}$ are such that:

\begin{equation*}
\begin{cases}
\bar{\lambda}^{*}=[\bar{\lambda}^{*}_{1} \quad \bar{\lambda}^{*}_{1}]\\
\widetilde{\lambda}^{*}=[\widetilde{\lambda}^{*}_{1} \quad \widetilde{\lambda}^{*}_{1}]
\end{cases}
\end{equation*}

where $\bar{\lambda}^{*}_{1} \in \Re^{2 \times 1}$ and $\widetilde{\lambda}^{*}_{1} \in \Re^{2 \times 1}$. Then, $R$ is a matrix with rank $1$, so the condition $\lambda_{R}<1$, yields:

\begin{equation*}
\begin{cases}
\lambda_{R}=trace(R)<1\\
R=\frac{1}{2} \cdot [(\bar{\lambda}^{*}+\widetilde{\lambda}^{*})+(\bar{\lambda}^{*}+\widetilde{\lambda}^{*})']
\end{cases}
\end{equation*}

where the operator $trace(\cdot)$ is the sum of the diagonal elements. The condition for $\lambda_{R}<1$ leads:

\begin{equation}\label{Conditions Hopf II}
\lambda_{R}<1 \Leftrightarrow \frac{1}{2} \cdot (\bar{\lambda}_{11}^{*}+\bar{\lambda}_{22}^{*}+\widetilde{\lambda}_{11}^{*}+\widetilde{\lambda}_{22}^{*})<1
\end{equation}

where $\bar{\lambda}^{*}=
\begin{bmatrix}
\bar{\lambda}_{11}^{*} & \bar{\lambda}_{12}^{*}\\
\bar{\lambda}_{21}^{*} & \bar{\lambda}_{22}^{*}
\end{bmatrix}$ and $\widetilde{\lambda}^{*}=
\begin{bmatrix}
\widetilde{\lambda}_{11}^{*} & \widetilde{\lambda}_{12}^{*}\\
\widetilde{\lambda}_{21}^{*} & \widetilde{\lambda}_{22}^{*}
\end{bmatrix}$.

Finally, the set of conditions to satisfy arise from equations (\ref{Conditions Hopf I}) and (\ref{Conditions Hopf II}):

\begin{equation}\label{Conditions Hopf III}
\begin{cases}
\mid 3 \cdot x^{2} + y^{2} \mid \leq \bar{\lambda}_{11}^{*}-\widetilde{\lambda}_{11}^{*}\\
\mid x^{2} +3 \cdot y^{2} \mid \leq \bar{\lambda}_{22}^{*}-\widetilde{\lambda}_{22}^{*}\\
\mid 2 \cdot x \cdot y \mid \leq \bar{\lambda}_{12}^{*}-\widetilde{\lambda}_{12}^{*}\\
\mid (\alpha \cdot x -y ) + 3 \cdot x \cdot (x^{2}+y^{2}) \mid \leq min_{i=1,2} \quad \widetilde{\lambda}^{*}_{1i}\\
\mid  (\alpha \cdot y + x ) + 3 \cdot y \cdot (x^{2}+y^{2}) \mid \leq min_{i=1,2} \quad \widetilde{\lambda}^{*}_{2i}\\
%
-\frac{1}{\alpha} \cdot (\bar{\lambda}_{11}^{*}+\bar{\lambda}_{22}^{*}+\widetilde{\lambda}_{11}^{*}+\widetilde{\lambda}_{22}^{*})<1
\end{cases}
\end{equation}

Noticing that the condition on $\lambda_{R}$ is limiting the values for $\{\bar{\lambda}_{ij}^{*},\widetilde{\lambda}^{*}_{ij}\}$ and $\{ \bar{\lambda}_{12}^{*},\bar{\lambda}_{21}^{*},\widetilde{\lambda}^{*}_{12},\widetilde{\lambda}^{*}_{21} \}$ are not present there, then it is possible to choose any value for them. If in particular they are such that:

\begin{equation*}
\begin{cases}
\widetilde{\lambda}^{*}_{12} \rightarrow +\infty\\
\widetilde{\lambda}^{*}_{21} \rightarrow +\infty\\
\end{cases}
\end{equation*}

In this way, conditions in equation (\ref{Conditions Hopf III}), lead:

\begin{equation*}
\begin{cases}
\mid 3 \cdot x^{2} + y^{2} \mid \leq \bar{\lambda}_{11}^{*}-\widetilde{\lambda}_{11}^{*}\\
\mid x^{2} +3 \cdot y^{2} \mid \leq \bar{\lambda}_{22}^{*}-\widetilde{\lambda}_{22}^{*}\\
\mid (\alpha \cdot x -y ) + 3 \cdot x \cdot (x^{2}+y^{2}) \mid \leq \widetilde{\lambda}^{*}_{11}\\
\mid  (\alpha \cdot y + x ) + 3 \cdot y \cdot (x^{2}+y^{2}) \mid \leq \widetilde{\lambda}^{*}_{22}\\
%
-\frac{1}{\alpha} \cdot (\bar{\lambda}_{11}^{*}+\bar{\lambda}_{22}^{*}+\widetilde{\lambda}_{11}^{*}+\widetilde{\lambda}_{22}^{*})<1
\end{cases}
\end{equation*}

One possibility is to choose $2 \cdot \widetilde{\lambda}_{11}^{*}=\bar{\lambda}_{11}^{*}=-\frac{\alpha}{3},\quad 2 \cdot \widetilde{\lambda}_{22}^{*}=\bar{\lambda}_{22}^{*}=-\frac{\alpha}{3}$, leading a domain of attraction $x\in[-0.08,0.08],\quad y\in[-0.08,0.08]$.

Next section is providing a set of examples to show the applicability of both: Theorem \ref{Main Theorem} and Corollary \ref{Systematic Conditions for Attraction}, together with some discussion about the accuracy when the domain of attraction obtained has to be used as an estimation to the truly global one.

\section{Examples}\label{Numerical Examples}

This section is devoted to show the applicability of the main theorem (Theorem \ref{Main Theorem}), as well as the its corollary. They together allow, in some cases, a precise estimation of the domain of attraction while asymptotic stability is established.

\begin{example}[An example with finite domain of attraction]

This first example is taken from the book \cite{Bacciotti05}, pp. 53:

\begin{equation*}
\begin{cases}
\dot{x_{1}(t)}=-x_{1}+2\cdot x_{1}^{2}\cdot x_{2}\\
\dot{x_{2}(t)}=-x_{2}
\end{cases}
\end{equation*}

For this system is known that the domain of attraction is given by $x_{1} \cdot x_{2} <1$, then this is a boundary for the estimation with Theorem \ref{Main Theorem}. Calculating the Jacobian:

\begin{equation}\label{Jacobian Ex 1}
\frac{\partial{f(x)}}{\partial{x}}=
\begin{bmatrix}
-1+4 \cdot x_{1} \cdot x_{2} & 2\cdot x_{1}^{2}\\
0 & -1
\end{bmatrix}
\end{equation}

Then, calculating at the origin:

\begin{equation*}
\frac{\partial{f(x)}}{\partial{x}}=
\begin{bmatrix}
-1 & 0\\
0 & -1
\end{bmatrix}
\end{equation*}

This Jacobian is clearly Hurwitz with matrix $P$ as follows:

\begin{equation*}
P=\frac{1}{2}\cdot
\begin{bmatrix}
1 & 0\\
0 & 1
\end{bmatrix}
\end{equation*}

On the other hand, the conditions on Theorem \ref{Main Theorem} leads:

\begin{equation}\label{Conditions Ex1 I}
\begin{cases}
\mid \frac{\partial{f(x)}}{\partial{x}}-F \mid \leq \bar{\lambda}^{*}-\widetilde{\lambda}^{*} \Leftrightarrow
\left |
\begin{matrix}
4\cdot x_{1}\cdot x_{2} & 2 \cdot x_{1}^{2}\\
0 & 0
\end{matrix} \right | \leq \bar{\lambda}^{*}-\widetilde{\lambda}^{*}\\
\mid \frac{\partial{f(x)}}{\partial{x}}\cdot x \mid \leq min_{i=1,2} \quad \widetilde{\lambda}^{*}_{i} \Leftrightarrow
\left |
\begin{matrix}
-x_{1}+4\cdot x_{1}^{2}\cdot x_{2} +2 \cdot x_{1}^{2}\cdot x_{2}\\
-x_{2}
\end{matrix} \right | \leq min_{i=1,2} \quad \widetilde{\lambda}^{*}_{i}
\end{cases}
\end{equation}

Analyzing the condition for $\lambda_{R}$:

\begin{equation*}
R=\frac{1}{2} \cdot [(\bar{\lambda}^{*}+\widetilde{\lambda}^{*})+(\bar{\lambda}^{*}+\widetilde{\lambda}^{*})']
\end{equation*}

On the other hand, if $\bar{\lambda}^{*},\widetilde{\lambda}^{*}$ are such that:

\begin{equation*}
\begin{cases}
\bar{\lambda}^{*}=[\bar{\lambda}^{*}_{1} \quad \bar{\lambda}^{*}_{1}]\\
\widetilde{\lambda}^{*}=[\widetilde{\lambda}^{*}_{1} \quad \widetilde{\lambda}^{*}_{1}]
\end{cases}
\end{equation*}

where $\bar{\lambda}^{*}_{1} \in \Re^{2 \times 1}$ and $\widetilde{\lambda}^{*}_{1} \in \Re^{2 \times 1}$. Then, $R$ is a matrix with rank $1$, so:

\begin{equation*}
\lambda_{R}=trace(R)<1
\end{equation*}

In this way, the second condition to satisfy arises:

\begin{equation}\label{Conditions Ex1 II}
\lambda_{R}<1 \Leftrightarrow \frac{1}{2} \cdot (\bar{\lambda}_{11}^{*}+\bar{\lambda}_{22}^{*}+\widetilde{\lambda}_{11}^{*}+\widetilde{\lambda}_{22}^{*})<1
\end{equation}

where $\bar{\lambda}^{*}=
\begin{bmatrix}
\bar{\lambda}_{11}^{*} & \bar{\lambda}_{12}^{*}\\
\bar{\lambda}_{21}^{*} & \bar{\lambda}_{22}^{*}
\end{bmatrix}$ and $\widetilde{\lambda}^{*}=
\begin{bmatrix}
\widetilde{\lambda}_{11}^{*} & \widetilde{\lambda}_{12}^{*}\\
\widetilde{\lambda}_{21}^{*} & \widetilde{\lambda}_{22}^{*}
\end{bmatrix}$. In summary, from equations (\ref{Conditions Ex1 I}) and (\ref{Conditions Ex1 II}):

\begin{equation}\label{Conditions Ex1 III}
\begin{cases}
\mid 4\cdot x_{1}\cdot x_{2}\mid \leq \bar{\lambda}_{11}^{*}-\widetilde{\lambda}_{11}^{*}\\
2 \cdot x_{1}^{2} \leq \bar{\lambda}_{12}^{*}-\widetilde{\lambda}_{12}^{*}\\
\mid -x_{1}+4\cdot x_{1}^{2}\cdot x_{2} +2 \cdot x_{1}^{2}\cdot x_{2} \mid \leq min_{i=1,2} \quad \widetilde{\lambda}^{*}_{1i}\\
\mid x_{2} \mid \leq min_{i=1,2} \quad \widetilde{\lambda}^{*}_{2i}\\
\frac{1}{2} \cdot (\bar{\lambda}_{11}^{*}+\bar{\lambda}_{22}^{*}+\widetilde{\lambda}_{11}^{*}+\widetilde{\lambda}_{22}^{*})<1
\end{cases}
\end{equation}

If in particular, it is chosen:

\begin{equation*}
\begin{cases}
\widetilde{\lambda}^{*}_{12} \rightarrow +\infty\\
\widetilde{\lambda}^{*}_{21} \rightarrow +\infty\\
\end{cases}
\end{equation*}

Then, the conditions in equation (\ref{Conditions Ex1 III}), yields:

\begin{equation*}
\begin{cases}
4 \cdot \mid x_{1}\cdot x_{2}\mid \leq \bar{\lambda}_{11}^{*}-\widetilde{\lambda}_{11}^{*}\\
\mid x_{1} \mid \cdot \mid -1+4 \cdot x_{2} +2 \cdot x_{1} \cdot x_{2} \mid \leq \widetilde{\lambda}^{*}_{11}\\
\mid x_{2} \mid \leq \widetilde{\lambda}^{*}_{22}\\
\frac{1}{2} \cdot (\bar{\lambda}_{11}^{*}+\bar{\lambda}_{22}^{*}+\widetilde{\lambda}_{11}^{*}+\widetilde{\lambda}_{22}^{*})<1\\
\forall x_{1} \in \Re
\end{cases}
\end{equation*}

Finally choosing $\bar{\lambda}_{11}^{*}-\widetilde{\lambda}_{11}^{*}=1$ and $\widetilde{\lambda}_{11}^{*}=\frac{1}{4},\widetilde{\lambda}_{22}^{*}=\frac{1}{2}$, this leads:

\begin{equation*}
\begin{cases}
\mid x_{1} \mid \cdot \mid -1+4 \cdot x_{2} +2 \cdot x_{1} \cdot x_{2} \mid \leq \frac{1}{4}\\
\mid x_{2} \mid \leq \frac{1}{2}\\
\end{cases}
\end{equation*}

The plot for this region is depicted in Figure \ref{Example 1 Domain Theorem 5}.  On the other hand, using Corollary \ref{Systematic Conditions for Attraction}, the domain of attraction is presented in Figure \ref{Example 1}.

\begin{figure}
\centering
\includegraphics[width=5cm]{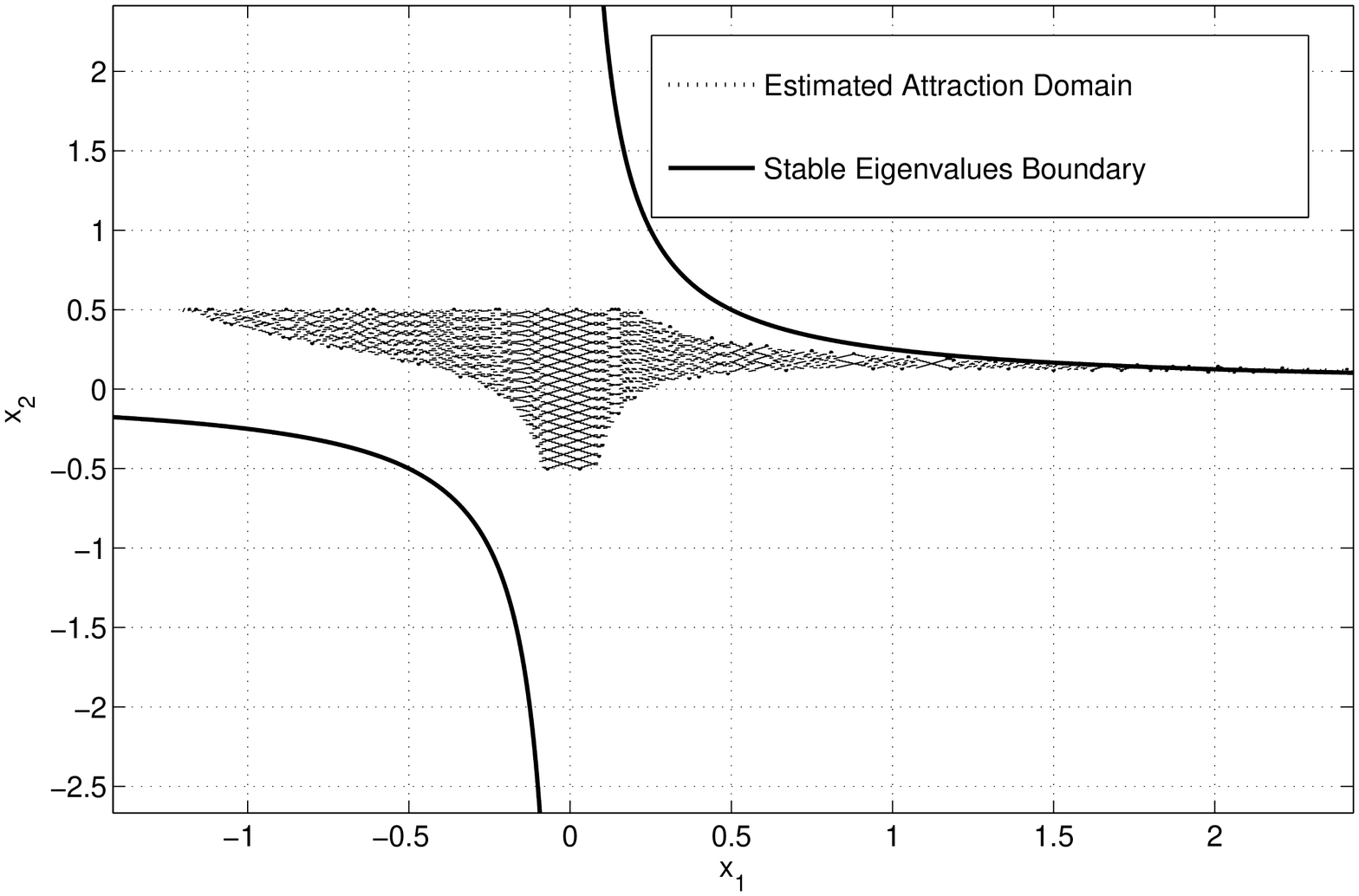}
\caption{Domain of attraction comparison using Theorem \ref{Main Theorem}.}\label{Example 1 Domain Theorem 5}
\end{figure}

\begin{figure}
\centering
\includegraphics[width=5cm]{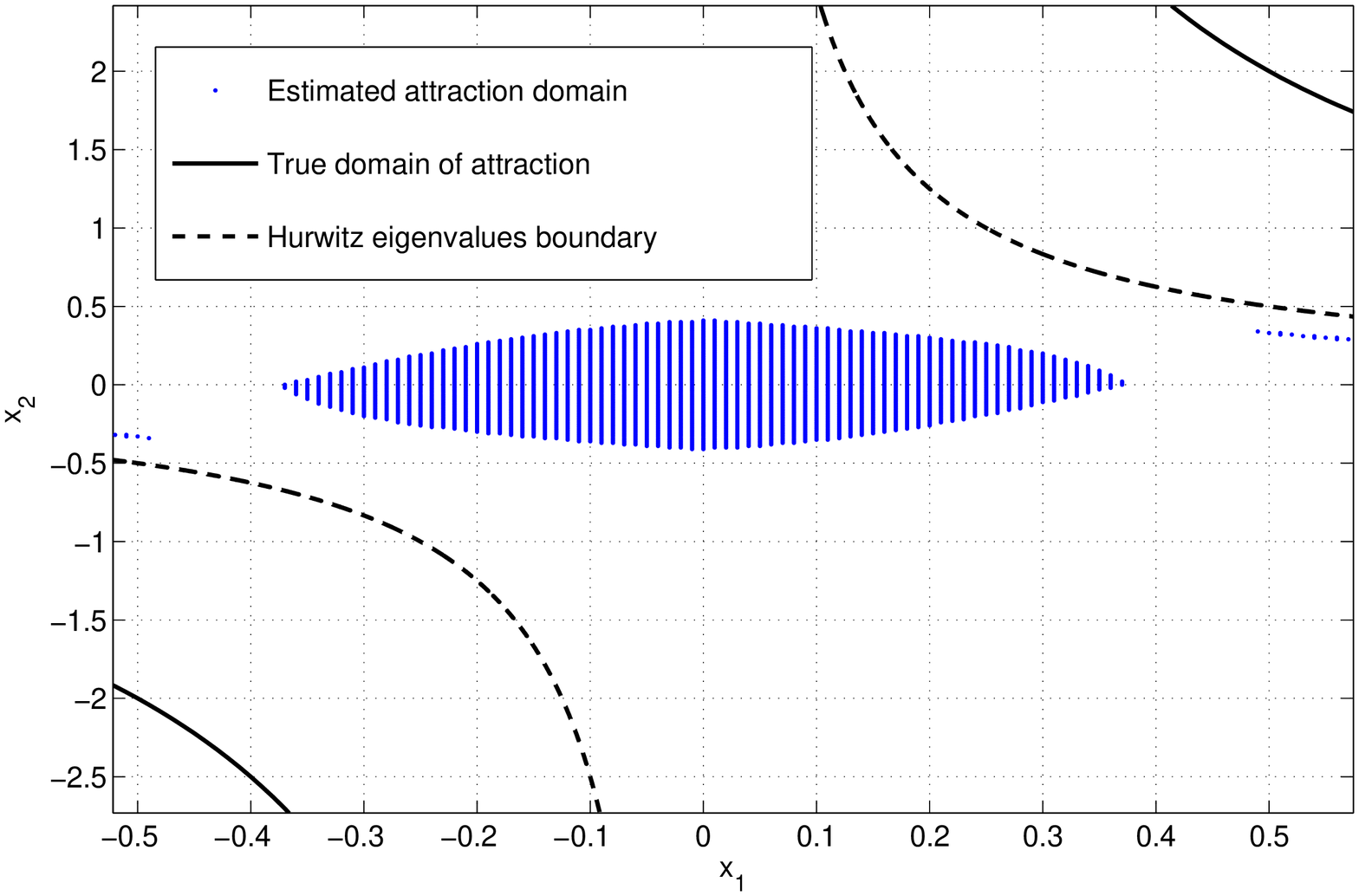}
\caption{Domain of Attraction Comparison using Corollary \ref{Systematic Conditions for Attraction}.}\label{Example 1}
\end{figure}

\end{example}

\begin{example}[The Van der Pol System]

The Van der Pol system is a very well known system which serves as a benchmark in many areas of science, in particular it is stable for negative values of the parameter $\mu$:

\begin{equation*}
\begin{cases}
\dot{x_{1}}(t)=x_{2}\\
\dot{x_{2}}(t)=-x_{1}+\mu \cdot x_{2} \cdot (1-x_{1}^{2})\\
\end{cases}
\end{equation*}

For this example, the systematic Corollary \ref{Systematic Conditions for Attraction} will be considered. To apply this result, it is necessary to obtain the jacobian:

\begin{equation*}
\frac{\partial{f(x)}}{\partial{x}}=
\begin{bmatrix}
0 & 1 \\
-1 +2 \mu \cdot x_{1} \cdot x_{2} & \mu \cdot (1-x_{1}^{2})
\end{bmatrix}
\end{equation*}

The eigenvalues are Hurwitz in the region depicted in Figure \ref{Example 2 Stable Eigenvalues}. It is clear that the stable region for the eigenvalues contains the equilibrium at the origin. In this way, the zone where the eigenvalue $\lambda_{R}$ in Corollary \ref{Systematic Conditions for Attraction} is less than one is plotted in Figure \ref{Example 2}, for $\mu=-1$.

\begin{figure}
\centering
\includegraphics[width=5cm]{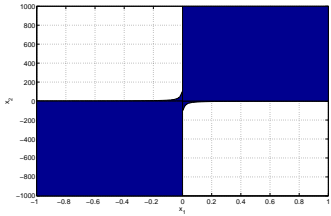}
\caption{Region for Hurwitz eigenvalues in the Van der Pol system.}\label{Example 2 Stable Eigenvalues}
\end{figure}

\begin{figure}
\centering
\includegraphics[width=5cm]{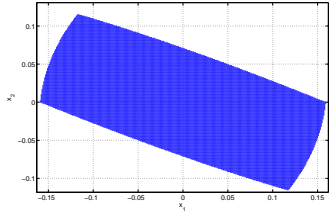}
\caption{Domain of Attraction obtained with Corollary \ref{Corollary of the Main Result}.}\label{Example 2}
\end{figure}

Since this system is \textit{stiff},\footnote{in the sense that possesses two main time-scales with different magnitude producing a sudden change into the direction motion of the system-this is one reason why it is used as a benchmark for numerical ODE solver testing} the obtained stability region is smaller than the whole domain of attraction depicted in Figure \ref{True Domain Attraction Ex 2}.

\begin{figure}
\centering
\includegraphics[width=5cm]{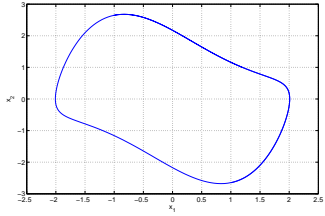}
\caption{True Domain of Attraction for the Stable Van der Pol system.}\label{True Domain Attraction Ex 2}
\end{figure}

\end{example}

\begin{example}[The Markus-Yamabe conjecture, discussion]

One special class of ODE's are those with Hurwitz jacobian matrix for all $\Re^{n}$. For these systems, the intuition could say they are always stable. However, as it is well known since 1995 this is not true for dimensions bigger than 2 in general cases(see \cite{Meisters96} and \cite{Chamberland04}).

Taking a look of Corollary \ref{Systematic Conditions for Attraction}, is clear that the region for stability is the intersection of two main regions:

\begin{itemize}
\item{The region where the eigenvalues of $\frac{\partial{f(x)}}{\partial{x}}$ are Hurwitz.}
\item{The region where $\lambda_{R}<1$.}
\end{itemize}

A conclusion using the theory in this paper is that the Markus-Yamabe conjecture was uncomplete, missing the second condition above.

\end{example}

\begin{example}[An example where the test is not applicable]

This example is due to Krasovskii (see \cite{Bacciotti05}, pp. 32, \cite{Krasowski55} and \cite{Auslander64}). This system has an equilibrium at the origin but the jacobian is unstable for $\Re^{2}$:

\begin{equation*}
\begin{cases}
\dot{x_{1}(t)}=x_{2}\\
\dot{x_{2}(t)}=-x_{1}+x_{2}\cdot (x_{1}^{2}+x_{2}^{2})^{3} \cdot sin(\frac{\pi}{(x_{1}^{2}+x_{2}^{2})})^{2}
\end{cases}
\end{equation*}

Since the jacobian is unstable, it is not possible to find an stable matrix $F$ as required by Theorem \ref{Main Theorem}, otherwise, in the view of Lemma \ref{Stability of A} the conclusion would be that the jacobian is stable leading a contradiction.

\end{example}

\begin{example}[Domain of attraction discussion]

Even when the focus on this paper is to derive an analytic tool to determine stability/attraction for Nonlinear ODE's, this example intends to discuss how precise is our prediction of the truly domain of attraction, in fact, the Normal-Hopf example in Section \ref{Sufficient Conditions for  ODE's} is very conservative but it is expected that the general scenario is of acceptable conservatism. In fact, it will be shown that for many systems the prediction exhibits high precision.

A class of systems with known domain of attraction are planar systems with circular unstable limit cycles. Those dynamics can be written in polar coordinates as follows:

\begin{equation*}
\begin{cases}
\dot{r(t)}=h(r)\\
\dot{\theta(t)}=1
\end{cases}
\end{equation*}

here $h(r)$ is negative for a finite interval on $r$ containing the origin ($r=0$), this idea is depicted in Figure \ref{Shape of h(r)}.

\begin{figure}
\centering
\includegraphics[width=5cm]{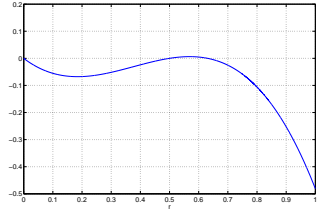}
\caption{Shape of $h(r)$}
\end{figure}\label{Shape of h(r)}

Notice that in order to have $f(x=0)=0$, it will be needed $h(0)=0$. Moreover for a finite circular domain of attraction, $h(\mu)=0$ with $\mu$ a finite positive value. In this way and calculating the Jacobian:

\begin{equation*}
\frac{\partial{f(x)}}{\partial{x}}=\frac{1}{r}\cdot
\begin{bmatrix}
r \cdot h(r) +cos(\theta)^{2}\cdot (h'(r)\cdot r-h(r))& - r + sin(\theta) \cdot cos(\theta) \cdot (h'(r)\cdot r-h(r))\\
r+ sin(\theta) \cdot cos(\theta) \cdot (h'(r)\cdot r-h(r)) & r \cdot h(r) +sin(\theta)^{2}\cdot (h'(r)\cdot r-h(r))
\end{bmatrix}
\end{equation*}

where $h'(r)=\frac{dh(r)}{dr}$. Since the border of the domain of attraction is circular, in order to analyze the condition in Corollary \ref{Systematic Conditions for Attraction}, $\lambda_{R}<1$, it can be settle:

\begin{equation*}
r=\mu,\quad \forall \theta \in \Re
\end{equation*}

In this way, the Jacobian for $\theta=0$\footnote{The rest of the cases lead equivalent results since is a circular domain.}:

\begin{equation*}
\frac{\partial{f(x)}}{\partial{x}}=
\begin{bmatrix}
h'(\mu) & -1\\
1 & 0
\end{bmatrix}
\end{equation*}

then:

\begin{equation*}
P=
\begin{bmatrix}
\frac{1}{h'(\mu)} & \frac{1}{2}\\
\frac{1}{2} & -\frac{2+h'(\mu)^{2}}{2 \cdot h'(\mu)}
\end{bmatrix}
\end{equation*}

Moreover:

\begin{equation*}
R=2 \cdot \mu \cdot (\mid P \mid \cdot \mid
\begin{bmatrix}
0 & 0\\
1 & 1
\end{bmatrix}+
\begin{bmatrix}
1 & 1\\
0 & 0
\end{bmatrix}\cdot \mid P \mid)
\end{equation*}

Finally:

\begin{equation*}
\lambda_{R}=2 \cdot \mu \cdot \left| \frac{h'(\mu)+2+h'(\mu)^{2}}{h'(\mu)} \right |
\end{equation*}

Calculating $\lambda_{R}<1$ it can be seen that for the interval $\mu \in [0, 0.2321]$, the estimation is perfect. Realize that this result is valid for any function $h(r)$ making valid the result for a wide range of systems.

\end{example}

\section{Conclusions}

Sufficient conditions for asymptotic stability and attraction for nonlinear ODE's were presented. These conditions are systematic to apply, moreover the main theorem provides a region where the ODE under analysis is asymptotically stable and a quadratic Lyapunov function. Realize that such a region serves as an estimation, sometimes coincident, of the domain of attraction and the Lyapunov function can extend this region by applying methods like the one in \cite{Burhcardt07}.

On the other hand, while Theorem \ref{Main Theorem} needs the input of a stable matrix $F$, Corollary \ref{Systematic Conditions for Attraction} utilizes the jacobian matrix of the given nonlinear ODE. The former result provides a subset of the domain of attraction but the estimation depends strongly on the choice of that matrix $F$. In this sense, the corollary was presented to overcome this inconvenient making the method systematic.

Several examples were presented for both Theorem \ref{Main Theorem} and Corollary \ref{Systematic Conditions for Attraction}, showing the applicability of such a results, moreover Example $5$ analyzed systems where the estimation of the domain of attraction is exact.

Many future directions can be depicted, in fact the eigenvalue $\lambda_{R}(x)$ in Corollary \ref{Systematic Conditions for Attraction} it seems to increase from $x=0$ until it reaches the value $1$ exactly when $x$ belong to the border of the subset of the domain of attraction obtained. If this is confirmed, then this quantity $\lambda_{R}$ serves to determine the borders of the estimation of the domain of attraction when is equal to one.

Besides, the theory supporting all the development in this paper is based on Theorem \ref{Stability Criterion} which in essence is a sufficient condition for attraction in CPWL ODE's. Notice that the theorem is proved utilizing a \textit{finite} improper integral, the finiteness is proved by requiring that all the matrices in each simplex of the CPWL under analysis being stable. This requirement could be in principle a bit conservative and should be improved allowing some of the matrices being unstable, then enlarging the attraction region provided by $\lambda_{R}$.

It will a topic for future research also the application of the methodologies on this paper to stabilize control systems also known as non-autonomous. In the case of control systems, one has to derive a controller but also to ensure that the system is stable at least for a finite region, in this sense, sufficient conditions for stability are the perfect scenario in order to obtain such a controllers.

\section{Appendix}

We present here the proof of those lemmas and theorems along the paper.

\begin{proof}[proof of Lemma \ref{Stability of A}]

Starting from the condition:

\begin{equation}\label{lambda aux 2}
\mid -F + A + B \mid \leq \overline{\lambda^{*}}\\
\end{equation}

In this way and considering a matrix $P'=[P_{1},P_{2},\ldots,P_{n}]'$ , $P_{j}=[P_{j1},P_{j2},\ldots,P_{jn}]$ and multiplying both sides of equation (\ref{lambda aux 2}) by $\mid P_{kj}\mid,\quad j,k=1,..,n$:

\begin{equation*}
\begin{cases}
\mid P_{kj} \mid  \cdot \mid -F_{jl} +A_{jl}+\bar{B}_{jl} \mid \leq \mid P_{kj} \mid \cdot \overline{\lambda}^{*}´_{jl}\\
j,l=1,..,n
\end{cases}
\end{equation*}

In matrix form:

\begin{equation}\label{P1}
\mid P \mid' \cdot \mid -F + A +B \mid \leq \mid P \mid' \cdot \bar{\lambda}^{*}\\
\end{equation}

In a similar way is also true:

\begin{equation}\label{P2}
\mid (-F+A^{(i)}+\bar{B}^{(i)})' \mid \cdot \mid P \mid \leq \overline{\lambda}^{*'} \cdot \mid P \mid\\
\end{equation}

Summing equations (\ref{P1}) and (\ref{P2}):

\begin{equation*}
\mid P \mid' \cdot \mid -F + A + B \mid +\mid -F + A + B' \mid \cdot  \mid P \mid \leq \overline{\lambda}^{*'} \cdot \mid P \mid+\mid P \mid' \cdot \overline{\lambda}^{*}\\
\end{equation*}

Realizing that $\mid P \cdot (-F+A+B) + (-F+A+B)'\cdot  \mid P \mid \leq \mid P \mid' \cdot \mid (-F+A+B) \mid +\mid (-F+A+B)'\mid \cdot  \mid P \mid$ and multiplying by $\mid x \mid,\quad x\in\Re^{n}$ at right and left:

\begin{equation*}
\mid x \mid' \cdot (\mid P \cdot (-F+A+B)+(-F+A+B)' \cdot  P \mid) \cdot \mid x \mid \leq \mid x \mid' \cdot (\overline{{\lambda}^{*}}' \cdot \mid P \mid+\mid P \mid' \cdot \overline{{\lambda}^{*}})\cdot \mid x \mid\\
\end{equation*}

Taking into account that $ x' \mid (\cdot) \mid x \leq \mid x \mid' \mid (\cdot) \mid \mid x \mid$, it is possible to write:

\begin{equation*}
\begin{cases}
x' \cdot \mid (P' \cdot (-F+A+B)+(-F+A+B)' \cdot  P) \mid \cdot x \leq \\
\mid x \mid' \cdot (\overline{{\lambda}^{*}}' \cdot \mid P \mid+\mid P \mid' \cdot \overline{{\lambda}^{*}})\cdot \mid x \mid\\
\end{cases}
\end{equation*}

Finally resorting to the property $\mid a + b\mid \geq \mid a \mid - \mid b \mid,\quad \forall a,b \in \Re$ and summing $x' \cdot (F' \cdot P+P' \cdot F) \cdot x$:

\begin{equation}\label{Useful for next Corollary}
\begin{cases}
x' \cdot (P' \cdot A + A \cdot  P) \cdot x < \\
\mid x \mid' \cdot (\overline{{\lambda}^{*}}' \cdot \mid P \mid+\mid P \mid' \cdot \overline{{\lambda}^{*}})\cdot \mid x \mid + x' \cdot (F' \cdot P+P' \cdot F) \cdot x)+x' \cdot \mid (P' \cdot B + B \cdot  P) \mid \cdot x\\
\end{cases}
\end{equation}

where the property $\mid a+b \mid \geq a+b,\forall a,b \in \Re$ was invoked. Since $F$ is stable if and only if there exists $P=P'>0$ such that $(P'\cdot F + F'\cdot P)=-I$, with $I$ the identity matrix (see \cite{Kwatny00}, pp. 22), also considering the definition of positive definiteness (see \cite{Bellman65}), in order to guarantee that the matrix $A$ is negative definite (or stable) is enough to require:

\begin{equation}\label{intermediate}
\mid x \mid' \cdot (\overline{{\lambda}^{*}}' \cdot \mid P \mid+\mid P \mid' \cdot \overline{{\lambda}^{*}})\cdot \mid x \mid + x' \cdot (F' \cdot P+P' \cdot F) \cdot x)+x' \cdot \mid P' \cdot B+B \cdot  P \mid \cdot x \leq 0
\end{equation}

Recalling that $F' \cdot P+P' \cdot F=-I$, then equation (\ref{intermediate}) looks like:

\begin{equation*}
\begin{cases}
\frac{\mid x \mid' \cdot (\overline{{\lambda}^{*}}' \cdot \mid P \mid+\mid P \mid' \cdot \overline{{\lambda}^{*}})\cdot \mid x \mid + x' \cdot \mid P' \cdot B + B \cdot  P \mid \cdot x }{x' \cdot x} \leq 1\\
x'\cdot x > 0
\end{cases}
\end{equation*}

Utilizing the condition:

\begin{equation}\label{lambda aux 3}
\mid B \mid \leq \widetilde{\lambda}^{*}
\end{equation}

The requirement for stability of $A$ will be:

\begin{equation}\label{intermidiate I}
\begin{cases}
\frac{\mid x \mid' \cdot \underbrace{((\overline{{\lambda}^{*}}+\widetilde{\lambda}^{*})' \cdot \mid P \mid+\mid P \mid' \cdot (\overline{{\lambda}^{*}}+\widetilde{\lambda}^{*}))}_{R}\cdot \mid x \mid}{x' \cdot x} \leq 1\\
x'\cdot x > 0
\end{cases}
\end{equation}

Since this inequality should be satisfied for all $x\in\Re^{n}$ in virtue of the definition of positiveness, the worst case leads:

\begin{equation}\label{Condition R}
max_{x} \frac{\mid x \mid' \cdot R \cdot \mid x \mid}{\mid x \mid ' \cdot \mid x \mid} \leq 1
\end{equation}

This is called Rayleigh's quotient and have important properties (see \cite{Bellman65}). To finish, realize that matrices $R$ and $\mid P \mid$ are positive (in the sense that all the elements are positive), so posses positive eigenvalues with positive corresponding eigenvector (see \cite{Bellman65}, pp. 72, Theorem 18).

In this way, the condition ensuring stability for $A$ leads:

\begin{equation*}
\lambda_{R}<1
\end{equation*}
This completes the proof.

\end{proof}

\begin{proof}[proof of Corollary \ref{Lyapunov functions for Linear Systems}]

Under the hypothesis of Lemma \ref{Stability of A}, and recalling the proof of this Lemma, in particular equation \ref{Useful for next Corollary} leads:

\begin{equation*}
x'\cdot(P \cdot A+A' \cdot P)\cdot x < 0
\end{equation*}

This completes the proof.

\end{proof}

\begin{proof}[proof of Theorem \ref{Main Theorem}]

The proof in what follows starts with the CPWL approximation bounds introduced in \cite{Julian99}:

\begin{equation*}
|f(x)-A^{(i)}\cdot x+B^{(i)}|\leq\lambda
\end{equation*}

where $A^{(i)} \cdot x+ B^{(i)}$ are the linear approximations for $f(x)$ in every simplex. This means that $A^{(i)}\cdot x+B^{(i)}$ is bounded by two vector fields: $\lambda+f(x)$ and $-\lambda+f(x)$, Figure \ref{Regiones de Ai y Bi} is depicting this idea. Notice that the dash region in that figure is the area containing the possible CPWL vector fields for a given $\lambda$.

\begin{figure}
\centering
\includegraphics[width=5cm]{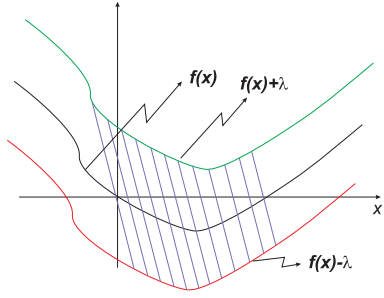}
\caption{Scheme of the Region Containing the CPWL approximation.}\label{Regiones de Ai y Bi}
\end{figure}

Since the focus is to create conditions over $f(x)$ related to $\{\tilde{\lambda^{*}},\bar{\lambda^{*}}\}$, we can consider the "worst case" $A^{(-\lambda)}\cdot x+B^{(-\lambda)}$ and $A^{(\lambda)}\cdot x+B^{(\lambda)}$ and show that those planes define a convex set \footnote{See \cite{Vanderbei99}, pp. 160-165} as follows (see Figure \ref{Planos de Peor Caso}):

\begin{equation}\label{Convex Set}
A^{(-\lambda)}\cdot x+B^{(-\lambda)} \leq f(x)-\lambda \leq A^{(i)}\cdot x+B^{(i)}\leq f(x)+\lambda \leq A^{( \lambda)}\cdot x+B^{\lambda)}
\end{equation}

The key idea it is to show that the inequalities above, are satisfied for at least one point inside each simplex and centering the attention in a generic simplex with vertices: $V_{i}\in \Re^{n},\quad i=1,\ldots,n+1$. In this way, if the gride size tends to zero\footnote{$V_{1}\rightarrow V_{2} and V_{i}\rightarrow V_{j},\quad i\neq j$}, it is possible to prove that: $A^{(\pm \lambda)}\rightarrow A^{(i)},\quad B^{(\pm \lambda)}\rightarrow B^{(i)}$.

Then, in place of using matrices $\{A^{(i)},B^{(i)}\}$, one uses the pair $\{A^{(\pm \lambda)},B^{(\pm \lambda)}\}$, which can be obtained for any generic nonlinear vector field $f(x)$ without the necessity of the evaluation of any CPWL approximation which moreover has to be reduced more and more in order to consider the limiting case when the grid size tend to zero.

Once the closed-form expressions for the pairs $\{A^{(\pm \lambda)},B^{(\pm \lambda)}\}$ have been obtained, it is possible to analyze attraction or stability as it going to be shown in what follows. For the sake of clarity we will only consider the analysis for $x \geq 0$, the case $x\leq 0$ shares conclusions.

\subsection*{$\Re^{1}$ analysis:}

In this way, let's calculate first these \textit{worst case} planes, then based on Figure \ref{Planos de Peor Caso}:

\begin{equation*}
\begin{cases}
P_{1}=f(V_{1})+\lambda,\quad P_{2}=f(V_{1})-\lambda,\quad P_{3}=f(V_{2})+\lambda,\quad P_{4}=f(V_{2})-\lambda\\
f(V_{1})+\lambda=A^{(\lambda)}\cdot V_{1}+B^{(\lambda)}\\
f(V_{2})+\lambda=A^{(\lambda)}\cdot V_{2}+B^{(\lambda)}\\
f(V_{1})-\lambda=A^{-(\lambda)}\cdot V_{1}+B^{-(\lambda)}\\
f(V_{2})-\lambda=A^{-(\lambda)}\cdot V_{2}+B^{-(\lambda)}
\end{cases}
\end{equation*}

In this way, the matrices $\{A^{(\pm \lambda)},B^{(\pm \lambda)}\}$ are given by:

\begin{equation}\label{A1, A2, B1, B2}
\begin{cases}
A^{(\lambda)}=\frac{f(V_{1})-f(V_{2})}{V_{1}-V_{2}}\\
B^{(\lambda)}=-[\frac{f(V_{1})-f(V_{2})}{V_{1}-V_{2}}]\cdot V_{1}+f(V_{1})+\lambda\\
A^{-(\lambda)}=\frac{f(V_{1})-f(V_{2})}{V_{1}-V_{2}}\\
B^{-(\lambda)}=-[\frac{f(V_{1})-f(V_{2})}{V_{1}-V_{2}}]\cdot V_{1}+f(V_{1})-\lambda\\
\end{cases}
\end{equation}

\begin{figure}
\centering
\includegraphics[width=5cm]{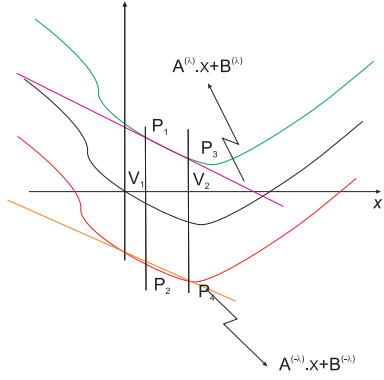}
\caption{"Worst case" planes.}\label{Planos de Peor Caso}
\end{figure}

At this point, the following conditions need to be verified:

\begin{equation}\label{Condition for a Convex Set}
\begin{cases}
A^{(-\lambda)}\cdot x+B^{(-\lambda)} \leq A^{(\lambda)}\cdot x+B^{(\lambda)}\\
A^{(-\lambda)}\cdot x+B^{(-\lambda)} \leq f(x)-\lambda\\
f(x)+\lambda \leq A^{(\lambda)}\cdot x+B^{(\lambda)}\\
\forall x \geq 0
\end{cases}
\end{equation}

Using the notation $V_{1}-V_{2}=\delta_{x}<0$ without loss of generality, the conditions above yields:

\begin{equation*}
\begin{cases}
A^{(-\lambda)}\cdot x+B^{(-\lambda)} \leq A^{(\lambda)}\cdot x+B^{(\lambda)} \Leftrightarrow 0 \leq 2 \cdot \lambda\\
A^{(-\lambda)}\cdot x+B^{(-\lambda)} \leq f(x)-\lambda \Leftrightarrow
f(x)-f(V_{1})\leq A^{(\lambda)}\cdot (x-V_{1}))\\
f(x)+\lambda \leq A^{(\lambda)}\cdot x+B^{(\lambda)}\Leftrightarrow
f(x)-f(V_{1})\geq A^{(-\lambda)}\cdot (x-V_{1})\\
\forall x \geq 0
\end{cases}
\end{equation*}

Considering equation (\ref{A1, A2, B1, B2}), it is clear that $A^{(\lambda)}=A^{(-\lambda)}$, then:

\begin{equation*}
\begin{cases}
0 \leq 2 \cdot \lambda\\
f(x)-f(V_{1})\leq A^{(\lambda)}\cdot (x-V_{1}))\\
f(x)-f(V_{1})\geq A^{(\lambda)}\cdot (x-V_{1})\\
\forall x \geq 0
\end{cases}
\end{equation*}

First condition is trivially satisfied taking into account that a CPWL approximation, always work with $\lambda \geq 0$, however, second condition means:

\begin{equation*}
\begin{cases}
f(x)-f(V_{1})=A^{(\lambda)}\cdot (x-V_{1}))\Leftrightarrow \frac{f(x)-f(V_{1})}{(x-V_{1}))}=A^{(\lambda)}\\
\forall x \geq 0
\end{cases}
\end{equation*}

Considering the limiting case $V_{1} \rightarrow V_{2}$, then:

\begin{equation*}
\begin{cases}
\frac{f(x)-f(V_{1})}{x-V_{1}}=lim_{V_{1} \rightarrow V_{2}}A^{(\lambda)}=\frac{\partial f(x)}{\partial x}|_{x=V_{1}}\\
\forall x \geq 0
\end{cases}
\end{equation*}

In the view of the \textit{mean value theorem}, this condition is always satisfied if $f(x)$ is continuous and with first derivative continuous - see \cite{Hubbard02}, pp. 148.

Once the validity of inequality (\ref{Convex Set})has been proved for at least one point $x\in [V_{1},V_{2}]$, it is missing to prove the following:

\begin{itemize}
\item{If $V_{1}\rightarrow V_{2}$, then: $A^{(\pm \lambda)}\rightarrow A^{(i)},\quad B^{(\pm \lambda)}\rightarrow B^{(i)}$.}
\item{The pair $\{A^{(\pm \lambda)},B^{(\pm \lambda)} \}$ satisfies: $\mid -F+A^{(\pm \lambda)}+B^{(\pm \lambda)} \mid \leq \bar{\lambda}^{*},\quad \mid B^{(\pm \lambda)} \mid \leq \widetilde{\lambda}^{*}$.}
\end{itemize}

The first item is going to be analyzed using inequality (\ref{Convex Set}) in the vertices $\{V_{1},V_{2}\}$ and recalling that $A^{(\lambda)}=A^{(-\lambda)}$:

\begin{equation*}
\begin{bmatrix}
V_{1} & 1\\
V_{2} & 1
\end{bmatrix}
\cdot
\begin{bmatrix}
A^{(\lambda)}\\
B^{(-\lambda)}
\end{bmatrix}
\leq
\begin{bmatrix}
V_{1} & 1\\
V_{2} & 1
\end{bmatrix}
\cdot
\begin{bmatrix}
A^{(i)}\\
B^{(i)}
\end{bmatrix}
\leq
\begin{bmatrix}
V_{1} & 1\\
V_{2} & 1
\end{bmatrix}
\cdot
\begin{bmatrix}
A^{(\lambda)}\\
B^{(\lambda)}
\end{bmatrix}
\end{equation*}

using absolute values:

\begin{equation*}
\left | T \cdot
\begin{bmatrix}
A^{(\lambda)}-A^{(i)}\\
h-B^{(i)}\\
\end{bmatrix} \right | \leq T \cdot
\begin{bmatrix}
0\\
\lambda
\end{bmatrix} \Leftrightarrow
|
\begin{bmatrix}
A^{(\lambda)}-A^{(i)}\\
h-B^{(i)}\\
\end{bmatrix} | \leq
\mid T^{-1} \mid \cdot T \cdot
\begin{bmatrix}
0\\
\lambda
\end{bmatrix}
\end{equation*}

where $T=\begin{bmatrix} V_{1} & 1\\ V_{2} & 1 \end{bmatrix}$ and $h=-A^{(\lambda)} \cdot V_{1}+f(V_{1})$. In this way, if $\lambda \rightarrow 0$, then:

\begin{equation*}
\left|\begin{bmatrix}
A^{(\lambda)}-A^{(i)}\\
h-B^{(i)}\\
\end{bmatrix}\right| =0
\end{equation*}
in other words: $A^{(\lambda)} \rightarrow A^{(i)}$ y $h \rightarrow B^{(i)}$.

\vspace{15mm}

The second item above, it going to be proved in what follows and states:

\begin{equation*}
\begin{cases}
\mid -F+A^{(\pm \lambda)}+B^{(\pm \lambda)} \mid \leq \bar{\lambda}^{*}\\
\mid B^{(\pm \lambda)} \mid \leq \widetilde{\lambda}^{*}\\
\end{cases}
\end{equation*}

For the sake of clarity in the exposition, it will be analyzed separately $B^{(\pm \lambda)}$ and $-F+A^{(\pm \lambda)}+B^{(\pm \lambda)}$:

\subsubsection*{For $B^{(\pm \lambda)}$:}\label{Analisis for Bi}

The requirement implies to satisfy the following:

\begin{equation*}
\begin{cases}
f(V_{1}) + \lambda - [\frac{f(V_{1})-f(V_{2})}{V_{1}-V_{2}}]\cdot V_{1} \leq \lambda_{2}^{*}\\
f(V_{1}) - \lambda - [\frac{f(V_{1})-f(V_{2})}{V_{1}-V_{2}}]\cdot V_{1} \geq -\lambda_{2}^{*}\\
\end{cases}
\end{equation*}

or, in other words:

\begin{equation}\label{Conditions for Bi}
\begin{cases}
f(V_{1}) \leq - \lambda + [\frac{f(V_{1})-f(V_{2})}{V_{1}-V_{2}}]\cdot V_{1} + \lambda_{2}^{*}\\
f(V_{1}) \geq  \lambda + [\frac{f(V_{1})-f(V_{2})}{V_{1}-V_{2}}]\cdot V_{1} -\lambda_{2}^{*}\\
\end{cases}
\end{equation}

\subsubsection*{For $-F+A^{(\pm \lambda)}+B^{(\pm \lambda)}$:}\label{Analisis for Ai}

This condition requires:

\begin{equation*}
\begin{cases}
-\lambda_{1}^{*} \leq -F+\frac{f(V_{1})-f(V_{2})}{V_{1}-V_{2}}\cdot (1-V_{1})+f(V_{1})-\lambda\\
-F+\frac{f(V_{1})-f(V_{2})}{V_{1}-V_{2}}\cdot (1-V_{1})+f(V_{1})+\lambda \leq \lambda_{1}^{*}\\
\end{cases}
\end{equation*}

This yields:

\begin{equation}\label{Conditions for Ai}
\begin{cases}
f(V_{1}) \leq \lambda_{1}^{*}+F-\frac{f(V_{1})-f(V_{2})}{V_{1}-V_{2}}\cdot (1-V_{1})-\lambda\\
f(V_{1}) \geq -\lambda_{1}^{*} +F-\frac{f(V_{1})-f(V_{2})}{V_{1}-V_{2}}\cdot (1-V_{1})+\lambda\\
\end{cases}
\end{equation}

Finally a region where the analysis is valid, has to be obtained.

\subsection*{Validity Regions in $\Re^{1}$:}

Equations (\ref{Conditions for Bi}) and (\ref{Conditions for Ai}) provide a region where the analysis is valid or in other words where there exists a CPWL approximation satisfying:

\begin{equation*}
\begin{cases}
\mid -F+A^{(\pm \lambda)}+B^{(\pm \lambda)} \mid \leq \bar{\lambda}^{*}\\
\mid B^{(\pm \lambda)} \mid \leq \widetilde{\lambda}^{*}\\
\end{cases}
\end{equation*}

In this way, the conclusions for stability and attraction valid for the CPWL vector field: $A^{(i)}\cdot x+B^{(i)}$ are also valid for $A^{(\pm \lambda)}\cdot x+B^{(\pm \lambda)}$ which can be immediately translated into conclusions for the nonlinear ODE.

For this reason, it is convenient to write together equations (\ref{Conditions for Bi}) and (\ref{Conditions for Ai}) to give:

\begin{equation}\label{Bounds Validity R1}
\begin{cases}
%
f(V_{1}) \leq - \lambda + [\frac{f(V_{1})-f(V_{2})}{V_{1}-V_{2}}]\cdot V_{1} + \lambda_{2}^{*}\\
f(V_{1}) \geq  \lambda + [\frac{f(V_{1})-f(V_{2})}{V_{1}-V_{2}}]\cdot V_{1} -\lambda_{2}^{*}\\
%
%
f(V_{1}) \leq \lambda_{1}^{*}+f(x=1)-\frac{f(V_{1})-f(V_{2})}{V_{1}-V_{2}}\cdot (1-V_{1})-\lambda\\
f(V_{1}) \geq -\lambda_{1}^{*} +f(x=1)-\frac{f(V_{1})-f(V_{2})}{V_{1}-V_{2}}\cdot (1-V_{1})+\lambda\\
\end{cases}
\end{equation}

The analysis will focus two main conditions:

\begin{itemize}\label{Requirements R1}
\item{The bounds for $f(V_{1})$ contain equilibrium points.}
\item{Suitable values for $\lambda$ in order to assure that a CPWL exists for our valid region in order to validate the proof using the CPWL theory theorem developed in Section \ref{Stability Criterion}}
\end{itemize}

The first item means a positive superior bound and a negative inferior one, performing this analysis to the conditions in (\ref{Bounds Validity R1}), it will be also obtained appropriate values for $\lambda$:

\begin{equation*}
\begin{cases}
%
\lambda \leq  \lambda_{2}^{*}+ [\frac{f(V_{1})-f(V_{2})}{V_{1}-V_{2}}]\cdot V_{1} \\
\lambda \leq  \lambda_{2}^{*}-[\frac{f(V_{1})-f(V_{2})}{V_{1}-V_{2}}]\cdot V_{1} \\
%
%
\lambda \leq \lambda_{1}^{*}-\frac{f(V_{1})-f(V_{2})}{V_{1}-V_{2}}\cdot (1-V_{1})+f(x=1)\\
\lambda \leq \lambda_{1}^{*}+\frac{f(V_{1})-f(V_{2})}{V_{1}-V_{2}}\cdot (1-V_{1})-f(x=1)\\
\end{cases}
\end{equation*}

Since $\lambda \geq 0$, considering the limiting case $V_{1} \rightarrow V_{2}$ and naming $V_{1}=x$, this leads:

\begin{equation*}
\begin{cases}
%
x \cdot \frac{\partial f(x)}{\partial{x}} \geq  -\lambda_{2}^{*}\\
x \cdot \frac{\partial f(x)}{\partial{x}} \leq  \lambda_{2}^{*}\\
%
%
(1-x)\cdot\frac{\partial f(x)}{\partial{x}} \leq \lambda_{1}^{*}+f(x=1)\\
(1-x)\cdot\frac{\partial f(x)}{\partial{x}} \geq -\lambda_{1}^{*}+f(x=1)\\
\end{cases}
\end{equation*}

In summary the region where the analysis is valid can be recast as:

\begin{equation*}
\begin{cases}
\mid x \cdot \frac{\partial f(x)}{\partial{x}} \mid \leq \lambda_{2}^{*}\\
\mid (1-x)\cdot\frac{\partial f(x)}{\partial{x}} -f(x=1)\mid \leq \lambda_{1}^{*}\\
\end{cases}
\end{equation*}

This completes the proof for the scalar case.

\subsection*{$\Re^{n}$ analysis:}

The development in this section extends the previous result for $\Re^{1}$ to vector fields in $\Re^{n}$, to do this, we consider $2\cdot n$ "worst case" planes defining a convex set as follows:

\begin{equation*}
\begin{cases}
A_{j}^{(-\lambda)} \cdot x +B_{j}^{(-\lambda)} \leq f_{j}(x)-\lambda_{j} \leq A^{(i)}_{j} \cdot x +B^{(i)}_{j} \leq f_{j}(x)+\lambda_{j} \leq A_{j}^{(\lambda)} \cdot x +B_{j}^{(\lambda)}\\
j=1,2,\ldots,n
\end{cases}
\end{equation*}

As for the $\Re^{1}$ case we need to show:

\begin{equation}\label{Condition for a Convex Set Rn}
\begin{cases}
A_{j}^{(-\lambda)}\cdot x+B_{j}^{(-\lambda)} \leq A_{j}^{(\lambda)}\cdot x+B_{j}^{(\lambda)}\\
A_{j}^{(-\lambda)}\cdot x+B_{j}^{(-\lambda)} \leq f_{j}(x)-\lambda\\
f_{j}(x)+\lambda \leq A_{j}^{(\lambda)}\cdot x+B_{j}^{(\lambda)}\\

j=1,2,\ldots,n\\
x=[x_{1},x_{2},\ldots,x_{n}]'\\
\forall x_{i} \geq 0
\end{cases}
\end{equation}

In this way, it is first needed $A_{j}^{(\pm \lambda)}$ and $B_{j}^{(\pm \lambda)}$ which can be obtained using $n+1$ points $\{ {x_{1},x_{2},\ldots,x_{n+1}} \}$ as follows:

\begin{equation*}
\begin{cases}
y_{k}=A_{j}^{(\pm \lambda)}\cdot x_{k}+B_{j}^{(\pm \lambda)}\\
A_{j}^{(\pm \lambda)}\in \Re^{1 \times n}\\
B_{j}^{(\pm \lambda)}\in \Re\\
j=1,\ldots,n\\
k=1,\ldots,n+1\\
\end{cases}\Leftrightarrow
\begin{bmatrix}
y_{1}\\
\vdots\\
y_{n+1}\\
\end{bmatrix}
=
\begin{bmatrix}
x_{1}'& 1\\
\vdots\\
x_{n+1}'& 1\\
\end{bmatrix}
\cdot
\begin{bmatrix}
A_{j}^{(\pm \lambda)'}\\
B_{j}^{(\pm \lambda)}\\
\end{bmatrix}
\end{equation*}

The set of points $\{ {y_{1},y_{2},\ldots,y_{n+1}} \}$ is related to $\{ {x_{1},x_{2},\ldots,x_{n+1}} \}$ by:

\begin{equation*}
\begin{cases}
y_{j}=f_{j}(V_{j})\pm \lambda_{j}\\
x_{j}=V_{j}\\
\lambda=[\lambda_{1},\lambda_{2},\ldots,\lambda_{n}]'\\
j=1,\ldots,n+1
\end{cases}
\end{equation*}

where $\{ {V_{1},V_{2},\ldots,V_{n+1}} \}$ is the set of vertices containing the general simplex $i^{th}$ considered for the analysis.

Then, it is possible to explicitly write $A^{\pm \lambda}_{j}$ and $B_{j}^{\pm \lambda}$ in matrix form:

\begin{equation}\label{Aj and Bj}
\begin{cases}
\begin{bmatrix}
A_{j}^{(\pm \lambda)'}\\
B_{j}^{(\pm \lambda)}\\
\end{bmatrix}
=
\begin{bmatrix}
V_{1}'& 1\\
\vdots\\
V_{n+1}'& 1\\
\end{bmatrix}^{-1}
\cdot
\begin{bmatrix}
f_{j}(V_{1}) \pm \lambda_{j}\\
\vdots\\
f_{j}(V_{n+1}) \pm \lambda_{j}
\end{bmatrix}\\
j=1,\ldots,n
\end{cases}
\end{equation}

Next, in order to effectively calculate $A_{j}^{(\pm \lambda)}$ and $B_{j}^{(\pm \lambda)}$ from equation (\ref{Aj and Bj}) we need the inverse of the matrix $\begin{bmatrix}
V_{1}'& 1\\
\vdots\\
V_{n+1}'& 1\\
\end{bmatrix}$.

A useful way of doing this is the one presented in \cite{Golubitsky99}, pp. 97, where applying row operations as in the well known method of Gauss for solving linear systems of algebraic equations it is possible to transform a matrix $M$ in $\widetilde{M}$ as follows:

\begin{equation*}
M=
\begin{bmatrix}
\begin{pmatrix}
V_{1}'& 1 \\
\vdots \\
V_{n+1}'& 1\\
\end{pmatrix}&
\begin{matrix}
|&\quad\\
|&I \\
|&\quad\\
\end{matrix}
\end{bmatrix}
\Longleftrightarrow
\widetilde{M}=
\begin{bmatrix}
\begin{matrix}
\quad & |\\
I & | \\
\quad & |\\
\end{matrix}
\begin{pmatrix}
V_{1}'& 1 \\
\vdots \\
V_{n+1}'& 1\\
\end{pmatrix}^{-1}&
\end{bmatrix}
\end{equation*}

where $I$ is the identity matrix. Then by subtracting the first with the second row in $M$, the second with the third and so on, we arrive to:

\begin{equation}\label{M}
M=
\begin{bmatrix}
\begin{pmatrix}
V_{1}'-V_{2}'& 0 \\
V_{2}'-V_{3}'& 0 \\
\vdots \\
V_{n}'-V_{n+1}'& 0 \\
V_{n+1}'& 1\\
\end{pmatrix}&
\begin{matrix}
\begin{matrix}
| & 1 &-1 & 0 & 0 & 0 & \ldots& 0\\
| & 0 &1 & -1 & 0 & 0 & \ldots& 0\\
| & 0 & 0 & 1 & -1 & 0 & \ldots& 0\\
| & \vdots\\
| & 0 & 0 & 0 & 0 & \ldots& 1 & -1\\
| & 0 & 0 & \quad & \quad & \ldots & 0 & 1
\end{matrix}
\end{matrix}
\end{bmatrix}
\end{equation}

On the other hand is possible to choose $n$ vertices forming a basis (see \cite{Julian99} for a further reading):

\begin{equation*}
V_{n+1}=
\begin{bmatrix}
V_{1} & V_{2} & \ldots & V_{n}
\end{bmatrix}
\cdot
\begin{bmatrix}
\gamma_{1}\\
\gamma_{2}\\
\vdots\\
\gamma_{n}
\end{bmatrix}
\end{equation*}

then there exists $\gamma^{*} \in \Re^{n \times 1}$ such that:

\begin{equation*}
-V_{n+1}=
\begin{matrix}
\begin{bmatrix}
V_{1}-V_{2} & V_{2}-V_{3} & \ldots & V_{n}-V_{n+1}
\end{bmatrix}
\cdot
\gamma^{*}
\end{matrix}
\end{equation*}

Now according the \textit{boundary configuration}\footnote{This is how the set of vertices is called in \cite{Julian99}} in \cite{Julian99}, pp. 58-65, equation (3.8) and more general in section (3.4.1), then it is possible to write every vertex for any simplex as follows:

\begin{equation*}
\begin{cases}
V_{i}=m_{i} \cdot \delta_{x_{i}} \cdot \breve{e_{i}}+\delta_{x_{i}} \cdot \sum^{i}_{j=1} (-1)^{c_{j}} \cdot \breve{e}_{j}\\
i=1,\ldots,n
\end{cases}
\end{equation*}

where $c_{j}\in{0,1}$, $m_{i}=1,2,\ldots,N_{i}$ with $N_{i}$ the number of simplicies used in each coordinate and $\delta_{x_{i}}$ their grid size . In this way:

\begin{equation*}
\begin{cases}
V_{i}-V_{i+1}=\delta_{x_{i}} \cdot (\sum^{i}_{j=1} (-1)^{c_{j}} \cdot \breve{e}_{j}- \sum^{i}_{j=1} (-1)^{c_{j}} \cdot \breve{e}_{j}-(-1) \cdot \breve{e}_{i+1}) \Rightarrow\\
V_{i}-V_{i+1}=\delta_{x_{i}} \cdot \breve{e}_{i+1}\\
i=1,\ldots,n
\end{cases}
\end{equation*}

providing the present analysis for $V_{i}>0, i=1,2,\ldots,n$-the rest of the cases are of similar consideration. This fact yields:

\begin{equation*}
Q=
\begin{matrix}
\begin{bmatrix}
V_{1}'-V_{2}'\\
\vdots\\
V_{n}'-V_{n+1}'
\end{bmatrix}
=
\begin{bmatrix}
\delta_{x_{1}} & 0 & \ldots & 0\\
0 & \delta_{x_{2}} & \ldots & 0\\
\vdots & \quad & \ddots \\
0 & 0 & \ldots & \delta_{x_{n}}
\end{bmatrix}
\end{matrix}
\end{equation*}

In this way the matrix $M$ from equation (\ref{M}) becomes:

\begin{equation*}
M=
\begin{bmatrix}
\begin{pmatrix}
\begin{matrix}
Q'
\end{matrix}
&
\begin{matrix}
0\\
\vdots\\
0
\end{matrix}\\
\begin{matrix}\\
\gamma^{*'}\cdot Q'
\end{matrix}
&
\begin{matrix}
\quad\\
1
\end{matrix}
\end{pmatrix}&
\begin{matrix}
\begin{matrix}
| & 1 &-1 & 0 & 0 & 0 & \ldots& 0\\
| & 0 &1 & -1 & 0 & 0 & \ldots& 0\\
| & 0 & 0 & 1 & -1 & 0 & \ldots& 0\\
| & \vdots\\
| & 0 & 0 & 0 & 0 & \ldots& 1 & -1\\
| & 0 & 0 & \quad & \quad & \ldots & 0 & 1
\end{matrix}
\end{matrix}
\end{bmatrix}
\end{equation*}

Performing row operations we have:

\begin{equation*}
M=
\begin{bmatrix}
\begin{pmatrix}
\begin{matrix}
Q'
\end{matrix}
&
\begin{matrix}
0\\
\vdots\\
0
\end{matrix}\\
\begin{matrix}\\
0 & \ldots & 0
\end{matrix}
&
\begin{matrix}
\quad\\
1
\end{matrix}
\end{pmatrix}&


\begin{pmatrix}
\begin{matrix}
I^{*}
\end{matrix}
&
\begin{matrix}
0\\
\vdots\\
0\\
-1
\end{matrix}\\
\begin{matrix}\\
-\gamma^{*'} \cdot I^{*}
\end{matrix}
&
\begin{matrix}
\quad\\
1+\gamma^{*}_{n}
\end{matrix}
\end{pmatrix}

\end{bmatrix}
\end{equation*}

where:

\begin{equation*}
I^{*}=
\begin{bmatrix}
1 & -1 & 0 & \ldots & 0\\
0 & 1 & -1 & \ldots & 0\\
\vdots\\
0 & 0 & 0 & \ldots & 1\\
\end{bmatrix}
\end{equation*}

Realizing that $-\gamma^{*} \cdot I^{*}=-[\gamma^{*}_{1} \quad (-\gamma^{*}_{1}+\gamma^{*}_{2}) \quad (-\gamma^{*}_{2}+\gamma^{*}_{3}) \ldots (-\gamma^{*}_{n-1}+\gamma^{*}_{n})]$,we obtain the desired inverse:

\begin{equation*}
\begin{matrix}
    \begin{bmatrix}
    V_{1}'& 1\\
    \vdots\\
    V_{n+1}'& 1\\
    \end{bmatrix}^{-1}&
    =&
    \begin{bmatrix}
        \frac{1}{\delta_{x_{1}}} &  0 & \ldots & 0 & 0\\
        0 &  \frac{1}{\delta_{x_{2}}} & \ldots & 0 & 0\\
        \vdots\\
        0 &  & \ldots & \frac{1}{\delta_{x_{n}}} & 0\\
        0 &  \ldots & 0 & 0 & 1
    \end{bmatrix}
    \cdot
    \begin{bmatrix}
        1 & -1 & 0 & \ldots & 0\\
        0 & 1 & -1 & \ldots & 0\\
        \vdots\\
        0 & 0 & \ldots & 1 & -1\\
        -\gamma^{*}_{1} & -(-\gamma^{*}_{1}+\gamma^{*}_{2}) & \ldots & -(-\gamma^{*}_{n-1}+\gamma^{*}_{n}) & (1+\gamma^{*}_{n})
    \end{bmatrix}
\end{matrix}
\end{equation*}

Finally, $A_{j}^{(\pm \lambda)}$ and $B_{j}^{(\pm \lambda)}$ arise:

\begin{equation}\label{Aj, Bj}
\begin{cases}
A_{j}^{(\pm \lambda)'}=
\underbrace{\begin{bmatrix}
\frac{f_{j}(V_{1})-f_{j}(V_{2})}{\delta_{x_{1}}}\\
\frac{f_{j}(V_{2})-f_{j}(V_{3})}{\delta_{x_{2}}}\\
\vdots\\
\frac{f_{j}(V_{n})-f_{j}(V_{n+1})}{\delta_{x_{n}}}
\end{bmatrix}}_{A^{*}_{j}}\\
B_{j}^{(\pm \lambda)}=[f_{j}(V_{n+1}) \pm \lambda_{j}]-\underbrace{\sum^{n}_{i=1} \gamma^{*}_{i} \cdot [f_{j}(V_{i})-f_{j}(V_{i+1})]}_{V_{n+1}'\cdot A_{j}^{(\pm \lambda)}}\\
\end{cases}
\end{equation}

where $\gamma^{*}_{i}= \frac {V_{(n+1)i}}{\delta_{x_{i}}}$ with $V_{(n+1)i}$ each component of the vertex $n+1$. Once $A_{j}^{(\pm \lambda)}$ and $B_{j}^{(\pm \lambda)}$ are already determined, the conditions for a convex set depicted in equation (\ref{Condition for a Convex Set Rn}) lead:

\begin{equation*}
\begin{cases}
A^{(-\lambda)}_{j}\cdot x+B^{(-\lambda)}_{j} \leq A^{(\lambda)}_{j}\cdot x+B^{(\lambda)}_{j} \Leftrightarrow
0 \leq (B^{(\lambda)}_{j}-B^{-(\lambda)}_{j}) \Leftrightarrow\\
\Leftrightarrow 2 \cdot \lambda_{j} \geq 0\\
f_{j}(x)-f_{j}(V_{n+1})=(x-V_{n+1})'\cdot A_{j}^{*}\\
j=1,2,\ldots,n
\end{cases}
\end{equation*}

As for $\Re^{1}$, considering the limiting case $V_{i} \rightarrow V_{i+1}$, then:

\begin{equation*}
\begin{cases}
A^{(-\lambda)}_{j}\cdot x+B^{(-\lambda)}_{j} \leq A^{(\lambda)}_{j}\cdot x+B^{(\lambda)}_{j} \Leftrightarrow
0 \leq (B^{(\lambda)}_{j}-B^{-(\lambda)}_{j}) \Leftrightarrow\\
\Leftrightarrow 2 \cdot \lambda_{j} \geq 0\\
f_{j}(x)-f_{j}(V_{n+1})=(x-V_{n+1})'\cdot \frac{\partial f_{j}(x)}{\partial x}|_{x=V_{n+1}}\\
j=1,2,\ldots,n
\end{cases}
\end{equation*}

The first condition above is always satisfied since all the CPWL approximations used along the proof provides a positive quantity $\lambda_{j}$, on the other hand the last requiremnet is also true for continuous vector fields $f(x)$ and with continuous first derivative in the view of the mean value theorem for vector fields in $\Re^{n}$ (see \cite{Hubbard02}, pp. 148).


On the other hand as for $\Re^{1}$, the following it is remaining to be proved:

\begin{itemize}
\item{If $V_{1} V_{2},\quad i=1,\ldots,n$, then: $A^{(\pm \lambda)} \rightarrow A^{(i)},\quad B^{(\pm \lambda)} \rightarrow B^{(i)}$.}
\item{$\{A^{(\pm \lambda)} ,B^{(\pm \lambda)} \}$ satisfy: $\mid F_{jk}-A_{jk}^{(.\pm \lambda)}-B_{j}^{(\pm \lambda)}\mid \leq \bar{\lambda}_{jk}^{*}\quad
\mid B^{(i)}_{j} \mid \leq \widetilde{{\lambda}}_{jk}^{*},\quad j,k=1,\ldots,n$.}
\end{itemize}

The first item can be proved utilizing the inequality (\ref{Condition for a Convex Set Rn}) in the $n+1$ vertices: $V_{i},\quad i=1,\ldots,n+1$ and noticing that $A^{(\lambda)}=A^{(-\lambda)}$ (see (\ref{Aj, Bj})):

\begin{equation*}
\begin{bmatrix}
V_{1}' & 1\\
V_{2}' & 1\\
\vdots\\
V_{n+1}' & 1
\end{bmatrix}
\cdot
\begin{bmatrix}
A^{(\lambda)'}_{j}\\
B^{(-\lambda)}_{j}
\end{bmatrix}
\leq
\begin{bmatrix}
V_{1}' & 1\\
V_{2}' & 1\\
\vdots\\
V_{n+1}' & 1
\end{bmatrix}
\cdot
\begin{bmatrix}
A^{(i)'}_{j}\\
B^{(i)}_{j}
\end{bmatrix}
\leq
\begin{bmatrix}
V_{1}' & 1\\
V_{2}' & 1\\
\vdots\\
V_{n+1}' & 1
\end{bmatrix}
\cdot
\begin{bmatrix}
A^{(\lambda)'}_{j}\\
B^{(\lambda)}_{j}
\end{bmatrix}
\end{equation*}

using absolute values:

\begin{equation*}
\left| T \cdot
\begin{bmatrix}
A^{(\lambda)'}_{j}-A^{(i)}_{j}\\
h_{j}-B^{(i)}_{j}\\
\end{bmatrix} \right | \leq T \cdot
\begin{bmatrix}
0\\
0\\
\vdots\\
\lambda_{j}
\end{bmatrix} \Leftrightarrow
\left |
\begin{bmatrix}
A^{(\lambda)'}_{j}-A^{(i)'}_{j}\\
h_{j}-B^{(i)}_{j}\\
\end{bmatrix} \right | \leq
\left | T^{-1} \right | \cdot T \cdot
\begin{bmatrix}
0\\
0\\
\vdots\\
0\\
\lambda_{j}
\end{bmatrix}
\end{equation*}

where:

\begin{equation*}
\begin{cases}
T=
\begin{bmatrix}
V_{1} & V_{2} & \ldots & V_{n+1}\\
1 & 1 & \ldots & 1
\end{bmatrix}'\\
h_{j}=-V_{n+1}' \cdot A^{(\lambda)}_{j} + f_{j}(V_{n+1})
\end{cases}
\end{equation*}

In this way, if $\lambda \rightarrow 0$, then:

\begin{equation*}
\left |
\begin{bmatrix}
A^{(\lambda)}_{j}-A^{(i)}_{j}\\
h_{j}-B^{(i)}_{j}\\
\end{bmatrix} \right |=0
\end{equation*}

in other words: $A^{(\lambda)} \rightarrow A^{(i)}$ y $h \rightarrow B^{(i)}$.

The second item yields:

\begin{equation*}
\begin{cases}
\mid F_{jk}-A_{jk}^{(i)}-B_{j}^{(i)}\mid \leq \bar{\lambda}_{jk}^{*}\\
\mid B^{(i)}_{j} \mid \leq \widetilde{}{\lambda}_{jk}^{*}\\
j,k=1,\ldots,n
\end{cases}
\end{equation*}

As for $\Re^{1}$, the analysis will be split into $B^{(\pm \lambda)}_{j}$ and $F_{jk}-A_{jk}^{(\pm \lambda)}-B_{j}^{(\pm \lambda)}$:

\subsubsection*{\textit{For $B^{(\pm \lambda)}_{j}$}}\label{Analisis for Bi Rn}

In this case, the following should be satisfied:

\begin{equation*}
\begin{cases}
f_{j}(V_{n+1}) + \lambda_{j} -  V_{n+1}'\cdot A_{j}^{*} \leq \widetilde{\lambda}_{j}^{*}\\
f_{j}(V_{n+1}) - \lambda_{j} -  V_{n+1}'\cdot A_{j}^{*} \geq -\widetilde{\lambda}_{j}^{*}\\
\end{cases}
\end{equation*}

Realize that $V_{n+1}'\cdot A_{j}^{*}$ play the role of $\frac{f(V_{1})-f(V_{2})}{\delta_{x}}$ for the $\Re^{1}$ case.

In this way:

\begin{equation}\label{Validity Region Bi Rn}
\begin{cases}
f_{j}(V_{n+1}) \leq - \lambda_{j} +  V_{n+1}'\cdot A_{j}^{*} + \widetilde{\lambda}_{j}^{*}\\
f_{j}(V_{n+1}) \geq \lambda_{j} +  V_{n+1}'\cdot A_{j}^{*}  -\widetilde{\lambda}_{j}^{*}\\
\end{cases}
\end{equation}

\subsubsection*{\textit{For $f_{j}(\breve{e}_{k})-A_{jk}^{(\pm \lambda)}-B_{j}^{(\pm \lambda)}$}}\label{Analisis for Ai Rn}

It has to be guaranteed:

\begin{equation*}
\begin{cases}
-f_{j}(\breve{e_{k}})+A_{jk}^{(\lambda)} + B_{j}^{(\lambda)} \leq \bar{\lambda}_{jk}^{*}\\
-f_{j}(\breve{e_{k}})+A_{jk}^{(-\lambda)} + B_{j}^{(-\lambda)} \geq -\bar{\lambda}_{jk}^{*}\\
k=1,2,\ldots,n
\end{cases}
\end{equation*}

Recalling equation (\ref{Aj, Bj}):

\begin{equation*}
\begin{cases}
A_{jk}^{*} +f_{j}(V_{n+1})+\lambda_{j}-V_{n+1}'\cdot A_{j}^{*} - f_{j}(\breve{e_{k}}) \leq \bar{\lambda}_{jk}^{*}\\
A_{jk}^{*} +f_{j}(V_{n+1}) - \lambda_{j}-V_{n+1}'\cdot A_{j}^{*} - f_{j}(\breve{e_{k}}) \geq -\bar{\lambda}_{jk}^{*}\\
k=1,2,\ldots,n
\end{cases}
\end{equation*}

this yields:

\begin{equation}\label{Validity Region Ai Rn}
\begin{cases}
f_{j}(V_{n+1}) \leq \bar{\lambda}_{jk}^{*} + f_{j}(\breve{e_{k}})- \lambda_{j}-A_{jk}^{*} +V_{n+1}'\cdot A_{j}^{*}\\
f_{j}(V_{n+1}) \geq -\bar{\lambda}_{jk}^{*} + f_{j}(\breve{e_{k}})+ \lambda_{j}-A_{jk}^{*} +V_{n+1}'\cdot A_{j}^{*}\\
k=1,2,\ldots,n
\end{cases}
\end{equation}

\subsection*{Validity Regions in $\Re^{n}$}

Comparing equations (\ref{Validity Region Bi Rn}) with (\ref{Conditions for Bi}) and (\ref{Validity Region Ai Rn}) with (\ref{Conditions for Ai}) for $\Re^{n}$ and $\Re^{1}$ respectively, we immediately see that the conclusions for $\Re^{1}$ can be readily extended to the present $\Re^{n}$ analysis to give:

\begin{equation}\label{Bounds Validity Rn}
\begin{cases}
%
V_{(n+1)i} \geq 0 \\
f_{j}(V_{n+1}) \leq - \lambda_{j} +  V_{n+1}'\cdot A_{j}^{*} + \widetilde{\lambda}_{jk}^{*}\\
f_{j}(V_{n+1}) \geq \lambda_{j} +  V_{n+1}'\cdot A_{j}^{*}  -\widetilde{\lambda}_{jk}^{*}\\
%
f_{j}(V_{n+1}) \leq \bar{\lambda}_{jk}^{*} + f_{j}(\breve{e_{k}})- \lambda_{j}-A_{jk}^{*} +V_{n+1}'\cdot A_{j}^{*}\\
f_{j}(V_{n+1}) \geq -\bar{\lambda}_{jk}^{*} + f_{j}(\breve{e_{k}})+ \lambda_{j}-A_{jk}^{*} +V_{n+1}'\cdot A_{j}^{*}\\
j,k=1,2,\ldots,n
\end{cases}
\end{equation}

As for $\Re^{1}$, the present $\Re^{n}$ needs to fit the requirement of the superior bound being positive and the inferior one negative in order to contain equilibriums of $f(x)$. In this way:

\begin{equation*}
\begin{cases}
\widetilde{\lambda}_{jk}^{*} -\lambda_{j} + V_{n+1}'\cdot A_{j}^{*} \geq 0\\
-\widetilde{\lambda}_{jk}^{*} +\lambda_{j} + V_{n+1}'\cdot A_{j}^{*} \leq 0\\
\bar{\lambda}_{jk}^{*}+f_{j}(\breve{e_{k}})-\lambda_{j}- A_{jk}^{*} + V_{n+1}'\cdot A_{j}^{*} \geq 0\\
-\bar{\lambda}_{jk}^{*}+f_{j}(\breve{e_{k}})+\lambda_{j}- A_{jk}^{*} + V_{n+1}'\cdot A_{j}^{*} \leq 0\\
j,k=1,2,\ldots,n
\end{cases}
\end{equation*}

The appropriate values for $\lambda_{j}$ arises:

\begin{equation*}
\begin{cases}
\lambda_{j} \leq \widetilde{\lambda}_{jk}^{*} + V_{n+1}'\cdot A_{j}^{*} \\
\lambda_{j} \leq \widetilde{\lambda}_{jk}^{*} - V_{n+1}'\cdot A_{j}^{*} \\
\lambda_{j} \leq \bar{\lambda}_{jk}^{*} + f_{j}(\breve{e_{k}})- A_{jk}^{*} + V_{n+1}'\cdot A_{j}^{*} \\
\lambda_{j} \leq \bar{\lambda}_{jk}^{*} - f_{j}(\breve{e_{k}})+ A_{jk}^{*} - V_{n+1}'\cdot A_{j}^{*}\\
j,k=1,2,\ldots,n
\end{cases}
\end{equation*}

Since $\lambda_{j} \geq 0$, writing in matrix form, considering the limiting case $V_{i} \rightarrow V_{i+1}$ and naming in general $V_{n+1}=x$, it is obtained:

\begin{equation*}
\begin{cases}
x^{*'}\cdot \frac{\partial{f_{j}(x)}}{\partial{x}} \geq -\widetilde{\lambda}_{j}^{*}\\
x^{*'}\cdot \frac{\partial{f_{j}(x)}}{\partial{x}} \leq \widetilde{\lambda}_{j}^{*}\\
(-x^{*'}+I)\cdot \frac{\partial{f_{j}(x)}}{\partial{x}} \geq -\bar{\lambda}_{j}^{*}+F_{k}\\
(-x^{*'}+I)\cdot \frac{\partial{f_{j}(x)}}{\partial{x}} \leq \bar{\lambda}_{j}^{*}+F_{k}\\
j=1,2,\ldots,n
\end{cases}
\end{equation*}

where $x^{*'}=[x',x',\ldots,x']'$ and $F=[F_{1}\quad F_{2}\quad \ldots F_{n}]$. The conditions above can be compactly written as follows choosing the minimum bound for the first two inequalities above:

\begin{equation}\label{Validity Region Rn}
\begin{cases}
\mid x^{*'}\cdot \frac{\partial{f(x)}}{\partial{x}}' \mid \leq min_{j=1,2,\ldots,n} \quad \widetilde{\lambda}_{j}^{*}\\
\mid (-x^{*'}+I)\cdot \frac{\partial{f(x)}}{\partial{x}}' -F \mid \leq \overline{\lambda}^{*'}\\
\frac{\partial{f(x)}}{\partial{x}}'=
\begin{bmatrix}
\frac{\partial{f_{1}(x)}}{\partial{x_{1}}} & \frac{\partial{f_{2}(x)}}{\partial{x_{1}}} & \ldots & \frac{\partial{f_{n}(x)}}{\partial{x_{1}}}\\
\frac{\partial{f_{1}(x)}}{\partial{x_{2}}} & \frac{\partial{f_{2}(x)}}{\partial{x_{2}}} & \ldots & \frac{\partial{f_{n}(x)}}{\partial{x_{2}}}\\
\vdots\\
\frac{\partial{f_{1}(x)}}{\partial{x}} & \frac{\partial{f_{2}(x)}}{\partial{x}} & \ldots & \frac{\partial{f_{n}(x)}}{\partial{x}}\\
\end{bmatrix}
\end{cases}
\end{equation}

These conditions can be further simplified as follows:

\begin{equation}\label{Validity Region Rn}
\begin{cases}
\mid x^{*'}\cdot \frac{\partial{f(x)}}{\partial{x}}' \mid \leq min_{j=1,2,\ldots,n} \quad \widetilde{\lambda}_{j}^{*}\\
\mid \frac{\partial{f(x)}}{\partial{x}} -F' \mid \leq \overline{\lambda}^{*}-\widetilde{\lambda}^{*}\\
\end{cases}
\end{equation}

To complete the proof, it is remarkable that the conditions obtained in this proof are independent of the coordinate system chosen.

Up to this point, it has been proved that there exists a region given by (\ref{Validity Region Rn}) and a precision CPWL under which all the matrices $A^{(i)}$ are Hurwitz in each simplex, because of that and since the objective it is to prove that the aforementioned region it is also a region where the nonlinear ODE: $\dot{x}(t)f(x)$ is attractive to its equilibrium point, it will be necessary to prove that this region is an invariant set.

By virtue of inequality (\ref{Condition for a Convex Set Rn}), if the trajectories of: $\dot{x(t)}=A^{(i)} \cdot x + B^{(i)}$, are considered:

\begin{equation*}
\begin{cases}
A^{(\lambda)}\cdot x +B^{(- \lambda)} \leq \underbrace{A^{(i)}\cdot x + B^{(i)} }_{\dot{x}(t)}\leq  A^{(\lambda)}\cdot x +B^{(\lambda)} \\
\mid B^{(\pm \lambda)} \mid \leq min_{i=1,2,\ldots,n}\quad \widetilde{\lambda}_{i}^{*}\\
\forall x \in \Omega\\
\Omega=\{x: \mid \frac{\partial{f(x)}}{\partial{x}}\cdot x \mid \leq min_{i=1,2,\ldots,n}\quad \widetilde{\lambda}_{i}^{*},\quad  \mid \frac{\partial{f(x)}}{\partial{x}} -F' \mid \leq (\bar{\lambda}^{*}-\widetilde{\lambda}^{*}) \}\\
A^{(\lambda)}=
\begin{bmatrix}
\frac{f(V_{1})-f(V_{2})}{\delta_{x_{1}}}\\
\frac{f(V_{2})-f(V_{3})}{\delta_{x_{2}}}\\
\vdots\\
\frac{f(V_{n})-f(V_{n+1})}{\delta_{x_{n}}}
\end{bmatrix}
\end{cases}
\end{equation*}

In order to solve the differential inequality above, the technique in \cite{Garcia09} is going to be applied to give:

\begin{equation*}
\begin{cases}
min_{\{1,2\}} \quad \{x^{1*}(t),x^{2*}(t)\} \leq x(t) \leq max_{\{1,2\}} \quad \{x^{1*}(t),x^{2*}(t)\} \\
\dot{x^{1*}}(t)=A^{(\lambda)} \cdot x^{1*}(t) + min_{i=1,2,\ldots,n}\quad \widetilde{\lambda}_{i}^{*}\\
\dot{x^{2*}}(t)=A^{(\lambda)} \cdot x^{1*}(t) - min_{i=1,2,\ldots,n}\quad \widetilde{\lambda}_{i}^{*}
\end{cases}
\end{equation*}

It is possible then, to assure that the trajectories given by the inequality: $A^{(\lambda)}\cdot x +B^{(- \lambda)} \leq \underbrace{A^{(i)}\cdot x + B^{(i)} }_{\dot{x}(t)}\leq  A^{(\lambda)}\cdot x +B^{(\lambda)}$ are contained in $\Omega,\quad \forall t \geq 0$ if it is ensured that the bounds are also contained there.

To do this, the extremums of $x^{1,2*}(t)$ are looked for, that means, the instants of time for which $x^{1,2*}(t)$ reach their maximums and minimums. In order to find these maximums or minimums, the elemental calculus result for which the instant of time where the extremum occur are utilized, in this way, it is necessary to evaluate $\dot{x}^{1,2*}(t)=0$ and $\dot{x}^{1,2*}(t)$ does not exist.

It is clear that the first possibility only happens if $t\rightarrow \infty$ in the view that $A^{(\lambda)}$ is Hurwitz, while the second one indicates $t=0$ due to the discontinuity introduced in considering only positive evolutions in time: ($t \in \Re^{+}$).

These reasons and considering the initial condition $x^{1,2*}(0)$ always inside the set $\Omega$, it will be enough to prove that: $lim_{t \rightarrow \infty }\quad x^{1,2*}(t) \in \Omega$, that is:

\begin{equation*}
\begin{cases}
lim_{t \rightarrow \infty }\quad x^{1*}(t)=-(A^{(\lambda)})^{-1} \cdot min_{i=1,2,\ldots,n}\quad \widetilde{\lambda}_{i}^{*}\\
lim_{t \rightarrow \infty }\quad x^{2*}(t)=(A^{(\lambda)})^{-1} \cdot min_{i=1,2,\ldots,n}\quad \widetilde{\lambda}_{i}^{*}
\end{cases}
\end{equation*}

equivalently:

\begin{equation*}
\begin{cases}
lim_{t \rightarrow \infty }\quad A^{(\lambda)} \cdot x^{1*}(t)=- \widetilde{\lambda}^{*} \leq min_{i=1,2,\ldots,n}\quad \widetilde{\lambda}_{i}^{*} \\
lim_{t \rightarrow \infty }\quad A^{(\lambda)} x^{2*}(t)=\widetilde{\lambda}^{*} \leq min_{i=1,2,\ldots,n}\quad \widetilde{\lambda}_{i}^{*}
\end{cases}
\end{equation*}

In this way, if $\lambda \rightarrow 0$:

\begin{equation*}
\mid \frac{\partial{f(x)}}{\partial{x}}\cdot x \mid \leq min_{i=1,2,\ldots,n}\quad \widetilde{\lambda}_{i}^{*}
\end{equation*}

This completes the proof.

\end{proof}

\begin{proof}[proof of Lemma \ref{Corollary of the Main Result}]

As it was proved in Theorem \ref{Main Theorem} each matrix $A^{(i)}$ is Hurwitz. In this way and applying the Corollary 4.1 in \cite{Johanson04} it is immediately seen that $V(x)=x'\cdot P \cdot x$ is a common quadratic Lyapunov function of:

\begin{equation*}
\dot{x^{(i)}}(t)=A^{(i)}\cdot x^{(i)}+ B^{(i)},\quad \forall i=1,2,\ldots N
\end{equation*}

where $N$ is the number of total simplices used for the partition of the CPWL ODE's considered. If, on the other hand, $\frac{\partial{P}}{\partial{x}}=0$ and using the bounds introduced in \cite{Julian99}, it is possible to write:

\begin{equation*}
\mid f(x^{(i)}) -A^{(i)}\cdot x^{(i)}- B^{(i)}\mid \leq \lambda,\quad \forall i=1,2,\ldots N
\end{equation*}

Multiplying both right sides by $\mid x^{(i)'} \mid \cdot \mid P \mid$:

\begin{equation*}
\mid x^{(i)'} \mid \cdot \mid P \mid \cdot \mid f(x^{(i)}) -A^{(i)}\cdot x^{(i)}- B^{(i)}\mid \leq \mid x^{(i)'} \mid \cdot \mid P \mid \lambda,\quad \forall i=1,2,\ldots N
\end{equation*}

Since $ \mid x^{(i)'} \cdot (P \cdot f(x^{(i)}) -A^{(i)}\cdot x^{(i)}- B^{(i)})\mid \leq \mid x^{(i)'} \mid \cdot \mid P \mid \cdot \mid f(x^{(i)}) -A^{(i)}\cdot x^{(i)}- B^{(i)}\mid$, it is obtained:

\begin{equation}\label{Useful Corollary I}
\mid x^{(i)'} \cdot P \cdot (f(x^{(i)}) -A^{(i)}\cdot x^{(i)}- B^{(i)})\mid  \leq \mid x^{(i)'} \mid \cdot \mid P \mid \lambda,\quad \forall i=1,2,\ldots N
\end{equation}

With a similar procedure is also valid:

\begin{equation}\label{Useful Corollary II}
\mid (f(x^{(i)}) -A^{(i)}\cdot x^{(i)}- B^{(i)})\cdot P \cdot x^{(i)}\mid  \leq \mid x^{(i)'} \mid \cdot \mid P \mid \lambda,\quad \forall i=1,2,\ldots N
\end{equation}

Finally summing equations (\ref{Useful Corollary I}) and (\ref{Useful Corollary II}):

\begin{equation*}
2 \cdot x^{(i)'} \cdot P \cdot f(x^{(i)}) -2 \cdot x^{(i)'} \cdot P \cdot(A^{(i)}\cdot x^{(i)}+ B^{(i)}) \leq 2 \cdot \mid x^{(i)'} \mid \cdot \mid P \mid \lambda,\quad \forall i=1,2,\ldots N
\end{equation*}

Since $V=x^{(i)'}\cdot P \cdot x^{(i)}$ is a Lyapunov function for $A^{(i)}\cdot x^{(i)}+ B^{(i)}$, then $2 \cdot x^{(i)'} \cdot P \cdot(A^{(i)}\cdot x^{(i)}+ B^{(i)})<0$. This is showing that:

\begin{equation*}
2 \cdot x^{(i)'} \cdot P \cdot f(x^{(i)}) -2 \cdot \mid x^{(i)'} \mid \cdot \mid P \mid \lambda \leq 2 \cdot x^{(i)'} \cdot P \cdot(A^{(i)}\cdot x^{(i)}+ B^{(i)}) < 0 ,\quad \forall i=1,2,\ldots N
\end{equation*}

If the limiting case $\lambda \rightarrow 0$ is considered, then in virtue of Error Bounds the Dynamics introduced in \cite{Garcia06}, it is known that the trajectories $x^{(i)}\rightarrow x(t)$, where $\dot{x(t)}=f(x)$, this leads:

\begin{equation*}
lim_{\lambda\rightarrow0}\quad 2 \cdot x^{(i)'} \cdot P \cdot f(x^{(i)}) -2 \cdot \mid x^{(i)'} \mid \cdot \mid P \mid \lambda < 0 \Leftrightarrow 2 \cdot x(t) \cdot P \cdot f(x(t)) < 0,\quad \forall i=1,2,\ldots N
\end{equation*}

This completes the proof.

\end{proof}

\begin{proof}[proof of Corollary \ref{Systematic Conditions for Attraction}]

As for the proof of Theorem \ref{Main Theorem}, it will be needed a Hurwitz matrix $F$ such that: $\mid \frac{\partial{f(x)}}{\partial{x}}-F\mid \leq \bar{\lambda}^{*}-\widetilde{\lambda}^{*}$. In this way if $\frac{\partial{f(x)}}{\partial{x}}$ is Hurwitz for all $x\in \Omega$ and using $\bar{\lambda}^{*}=\widetilde{\lambda}^{*}$, then if $\Omega$ is as follows:

\begin{equation*}
\begin{cases}
\Omega=\{ x: \lambda_{R} <1 \}\\
R = 2 \cdot (\underbrace{\mid P \mid \cdot \mid \frac{\partial{f(x)}}{\partial{x}} \cdot x \mid \cdot [1 \quad 1 \quad \ldots 1]}_{H(x)} + H(x)')\\
\text{with $\lambda_{R}$ the biggest eigenvalue of $R$}
\end{cases}
\end{equation*}

the ODE: $\dot{x}(t)=f(x)$ is attractive to the origin (equilibrium point) in the view of Theorem \ref{Main Theorem}.

This completes the proof.

\end{proof}

\section{Acknowledgments}

The authors would like to acknowledge Universidad Nacional del Sur, CIC, CONCIET and ANPCyT.


\bibliographystyle{plain}
\bibliography{references}

\end{document}